\documentclass[11pt,a4paper]{amsart}
\usepackage{fontenc,amsmath,mathrsfs,amsfonts,amsthm, verbatim,amssymb, MnSymbol, enumerate, multicol, hyperref, geroENG, color}
\usepackage[english]{babel}
\allowdisplaybreaks
\usepackage{times}
\usepackage[top=2.0cm, bottom=2.0cm, left=2.0cm, right=2.0cm]{geometry}

\begin{document}
\title{Continued fractions with large prime partial quotients}

\author[G. Gonz\'alez Robert]{Gerardo Gonz\'alez Robert}
\address{Department of Mathematical and Physical Sciences,  La Trobe University, Bendigo 3552, Australia. }
\email{G.Robert@latrobe.edu.au}

\author[M. Hussain]{Mumtaz Hussain}
\address{Department of Mathematical and Physical Sciences,  La Trobe University, Bendigo 3552, Australia. }
\email{m.hussain@latrobe.edu.au}

\author[B. Ward]{ Benjamin Ward}
\address{Benjamin Ward,  {Department of Mathematics, University of York, Heslington, YO10 5DD.} }
\email{benjamin.ward@york.ac.uk, ward.ben1994@gmail.com }

\author[L. White]{Lauren White}
\address{Department of Mathematical and Physical Sciences,  La Trobe University, Bendigo 3552, Australia. }
\email{Lauren.White@latrobe.edu.au}

\thanks{This research is supported by the ARC Discovery Project 200100994. The third named author is a Leverhulme Early Career Research Fellow.}

\begin{abstract} 
We determine the Lebesgue measure and Hausdorff dimension of various sets of real numbers with infinitely many partial quotients that are both large and prime, thus extending the well-known theorems by {\L}uczak (1997) and Huang-Wu-Xu (2020). 
To this end, we obtain new asymptotics on the tail end of the almost prime zeta function.
Our results include some recent work by Schindler-Zweim{\"u}ller (2023). 

\end{abstract}
\maketitle

\section{Introduction}

Continued fraction theory is a cornerstone in Diophantine approximation. Within the theory of Diophantine approximation, one of the first results we are likely to encounter is Dirichlet's Theorem.
\begin{teo01}[Dirichlet's Theorem]
Let $x \in \RE$.
For any $N\in\N$ there exists a pair $(p,q)\in \Za \times\N$ such that
\begin{equation*}
    \left|x-\tfrac{p}{q}\right|<\frac{1}{qN} \, ,\quad \text{ and } \, \quad 0<q\leq N\, .
\end{equation*}
\end{teo01}
An immediate implication is the following corollary.
\begin{coro01}[Dirichlet's Corollary]
Let $x \in \RE\backslash\Q$. Then there exists infinitely many pairs $(p,q)\in\Za\times\N$ such that 
\begin{equation*}
    \left|x-\frac{p}{q}\right|<\frac{1}{q^{2}}\, .
\end{equation*}
\end{coro01}

The results of Dirichlet are insightful in that they tell us that any irrational number has rational points within some prescribed neighbourhood. However, we are no closer to knowing what these rational numbers look like. It also does not tell us whether the constant of $1$ appearing on the right-hand side of the equations can be improved or not. The numbers for which Dirichlet's Corollary can be improved by an arbitrarily small constant are called \textit{well approximable}, and the numbers in the complementary set are called \textit{badly approximable}. The set of points for which the constant $1$ appearing in Dirichlet's theorem can be improved is called the set of \textit{Dirichlet improvable numbers}. Studying these sets is when the power of continued fraction theory kicks in. \par 

Take an irrational number $x\in(0,1)$ with continued fraction expansion $x=[a_1(x),a_2(x),\ldots]$ and let $(p_n)_{n\geq 1}$ and $(q_n)_{n\geq 1}$ be the sequence of numerators and denominators of its convergents.
The natural numbers $a_1(x),a_2(x), \ldots$ are the \textit{partial quotients} of $x$. For any $x\in(0,1)$, the partial quotients can be constructed via the Gauss map. It is well-known \cite[Theorems 9 \& 13]{KhinchinCFBook} that for all $n\in\N$ we have
\begin{equation*}
    \frac{1}{(a_{n+1}+2)q_{n}^{2}}<\left|x-\frac{p_{n}}{q_{n}}\right|<\frac{1}{a_{n+1}q_{n}^{2}}\, .
\end{equation*}
Hence, many approximation properties of $x$ are determined by its partial quotients. 
For example, if $(a_n(x))_{n\geq 1}$ is bounded, the inequalities above along with other basic approximation properties of the convergents imply that $x$ is badly approximable.  
Conversely, if the partial quotients are unbounded, then $x$ is well-approximable. Furthermore, it was proven by Davenport and Schmidt \cite{MR272722} that a number is Dirichlet improvable if and only if its partial quotients are bounded. Thus, the set of Dirichlet improvable numbers is precisely the set of badly approximable numbers and rational numbers.\par 

These relationships between the properties of partial quotients and sets of a Diophantine interest motivate the following setup. For function $\varphi:\N\to[1, \infty)$ define
\begin{equation*}
    \clE_{1}(\varphi):=\left\{ x \in [0,1]: a_{n}(x)\geq \varphi(n) \ \text{for i.m.}\,  \, n\in \N\right\}\, .
\end{equation*}
In the above definition and throughout, ``i.m.'' means ``infinitely many''. 
The set $\clE_{1}(\varphi)$ and its variations have been vastly studied, with one of the first results on its size due to Borel and Bernstein, see, for example, \cite[Theorem 1.1]{MR2136100}. 
Let $\leb$ denote the Lebesgue measure.
\begin{teo01}[Borel-Bernstein Theorem]\label{borel bernstein theorem}
Let $\varphi:\N\to [1, \infty)$. Then
\begin{equation*}
    \leb(\clE_{1}(\varphi))=\begin{cases}
        0 \quad \text{\rm if} \quad \displaystyle\sum_{n=1}^{\infty}\frac{1}{\varphi(n)}<\infty\, , \\
        1 \quad \text{\rm if} \quad \displaystyle\sum_{n=1}^{\infty}\frac{1}{\varphi(n)}=\infty\, . \\
    \end{cases}
\end{equation*}
\end{teo01}
The Borel-Bernstein Theorem applied to the function $\vphi(n)=n\log n$ tells us that the set of badly approximable numbers, and thus also the set of Dirichlet-improvable numbers by \cite{MR272722}, is Lebesgue-null. Furthermore, the complement (the set of well approximable numbers) has full Lebesgue measure. \par

Given a non-increasing function $\psi:\N\to \RE_{>0}$, Kleinbock and Wadleigh \cite{KleinbockWadleigh2018} defined the set of \textit{$\psi$-Dirichlet improvable numbers} as the set of those $x\in\RE$ such that for every large $N\in\N$ there is a pair $(p,q)\in\Za\times \N$ satisfying
\begin{equation*}
    \left|qx - p\right|< \psi(N) \quad \text{ and }\quad q\leq N \, .
\end{equation*}
Remarkably, this set can also be characterised by the properties of the partial quotients. Kleinbock and Wadleigh \cite[Lemma 2.2]{KleinbockWadleigh2018} showed that a number $x$ is not $\psi$-Dirichlet improvable if $a_n(x)a_{n+1}(x)\geq \tfrac{1}{q_n\psi(q_n)}-1$ for infinitely many $n\in\N$. 
This leads to the study of the metric properties of the set
\begin{equation*}
    \clE_{\ell}(\varphi):=\left\{ x \in [0,1]: a_{n}(x)\ldots a_{n+\ell-1}(x)\geq \varphi(n) \ \text{for i.m.}\,  \, n\in \N\right\}\, .
\end{equation*}
Kleinbock and Wadleigh gave a zero-one law for the Lebesgue measure of $\clE_2(\vphi)$ assuming $\vphi(n)\to \infty$ as $n\to\infty$.
As a consequence, they obtained the Lebesgue measure of the set of $\psi$-Dirichlet improvable numbers when $\psi$ is non-increasing and $q\psi(q)<1$ for large $q\in\N$. 
In \cite{HuangWuXu2020}, Huang, Wu, and Xu obtained $\leb(\clE_{\ell})$ for $\ell\in\N_{\geq 2}$ without the limit hypothesis on $\vphi$.

\begin{teo01}[{\cite[Theorem 1.5]{HuangWuXu2020}}]\label{huang wu xu lebesgue theory}
Let $\varphi:\N\to [2, \infty)$ be a positive function. Then
\begin{equation*}
    \leb(\clE_{\ell}(\varphi))=\begin{cases}
        0\quad \text{\rm if } \quad \displaystyle\sum_{n=1}^{\infty}\frac{\left(\log \varphi(n)\right)^{\ell-1}}{\varphi(n)}<\infty\, , \\
        1\quad \text{\rm if } \quad \displaystyle\sum_{n=1}^{\infty}\frac{\left(\log \varphi(n)\right)^{\ell-1}}{\varphi(n)}=\infty\, . \\
    \end{cases}
\end{equation*}
\end{teo01}
When we are in the convergence case, that is, $\clE_{\ell}(\vphi)$ is Lebesgue-null, we appeal to the finer notion of size, for instance, the Hausdorff dimension, denoted by $\hdim$ to distinguish between them.
One of the first results in this regard is due to {\L}uczak. 
Considering functions $\vphi_{b,c}:\N\to\RE_{>0}$ given by $\vphi_{b,c}(n)=c^{b^n}$, $n\in\N$, with $b,c>1$ fixed the following was proven.
\begin{teo01}[\cite{MR1464376, Luc1997}]\label{TEO:LuczakTheorem}
For every $c,b>1$, we have
\[
\hdim \left\{ x\in[0,1]: \forall n\in\N \quad a_{n}(x)\geq \varphi_{b,c}(n) \right\}
=
\hdim \clE_1(\vphi_{b,c}) 
=
\frac{1}{b+1} \, .
\]
\end{teo01}
The Hausdorff dimension theory for general functions $\vphi$ was obtained for the case $\ell=1$ by Wang and Wu \cite{WangWu2008} and it was later generalised to arbitrary $\ell\in \N$ by Huang, Wu, and Xu \cite{HuangWuXu2020}.
For any function $\vphi:\N \to\RE_{>0}$, define the non-negative quantities $B_{\vphi}$ and $b_{\vphi}$ by
\begin{equation*}
    \log B_{\vphi}:=\liminf_{n\to \infty}\frac{\log\varphi(n)}{n}\,\quad \text{ and } \quad \log b_{\vphi}:=\liminf_{n\to \infty} \frac{\log\log \varphi(n)}{n} \, . 
\end{equation*}
\begin{teo01}[{\cite[Theorem 1.7]{HuangWuXu2020}}]\label{TEO:HWX:HDIM}
Define the sequence of functions $(f_{\ell})_{\ell\geq 1}$ from $(0,1)$ to $\RE_{>0}$ by
\begin{equation*}
    f_{1}(s)=s\, \quad \text{ and } \quad f_{\ell}(s)=\frac{sf_{\ell-1}(s)}{1-s+f_{\ell-1}(s)} \, , \, \text{ for } \, \, \ell\geq 2\, .
\end{equation*}
Let $\vphi:\N\to \RE_{>0}$ be arbitrary. 
Then, for every $\ell\in\N$,
\begin{equation*}
    \hdim \clE_{\ell}(\varphi)
    =\begin{cases}
        1\, , \quad \hspace{7cm}\text{\rm if } \quad B_{\vphi}=1\, ,\\
        \inf\{s:P(T,-f_{\ell}(s)\log B_{\vphi}-s\log|T'|)\leq 0\}\, , \,\,\,\, \,  \quad \text{\rm if } \quad 1<B_{\vphi}<\infty \, ,\\
        \frac{1}{b_{\vphi}+1}\, , \quad \hspace{6.5cm}\text{\rm if } \quad B_{\vphi}=\infty\,.
    \end{cases}
\end{equation*}
\end{teo01}
The definition of the pressure function $P(T,\cdot)$ in Theorem \ref{TEO:HWX:HDIM} is given in Subsection \ref{Subsec:PresureFunctions}. Without going into technical details, in our setup, the pressure function gives the Hausdorff dimension number as a continuous function of $B_{\vphi}$. It is $1$ when $B_{\vphi}\to 1$ and $1/2$ when $B_{\vphi}\to\infty$, see Proposition \ref{LE:DIMENSIONALNUMBER}.  

The above results have ignited an upsurge of papers on the metric properties of various sets in the theory of continued fractions. For example, one can look at:
\begin{itemize}
    \item the size of infinitely many partial quotients growing at least or at most at a certain rate \cite{MR4824910},
    \item the size of products of  partial quotients from arithmetic progressions \cite{HussainShulgaPisa},
    \item the size of the partial quotients being bounded by some path-dependent function \cite{MR3492979},
    \item the product of consecutive partial quotients belongs to designated intervals \cite{MR4671519, HussainShulgaPisa, MR4577484, MR4687014}.
\end{itemize} 
Alternatively, one can impose global conditions on the partial quotients. For instance, one can:
\begin{itemize}
    \item restrict all the partial quotients to $\{1,2,\ldots,M\}$ for some $M\in N$ \cite{Hensley1992},
    \item restrict all the partial quotients to $\{M,M+1,\ldots\}$ for some $M\in\N$ \cite{Cusick1971},
    \item restrict all even partial quotients to $\{1,\ldots,M\}$ for some $M\in\N$ \cite{HanciTurek2023},
    \item restrict all partial quotients to arithmetic progressions, prime numbers, squares, powers of 2 or 3, and lacunary sequences, see \cite[Table 1]{MR4068257}, see also \cite{MR4900954}. 
\end{itemize}
In this article we consider in some sense a combination of the two settings, by studying the set of continued fractions with infinitely many large prime partial quotients.
Firstly, we define the following functions: for $n\in\N$, let $a_n':(0,1)\setminus\QU\to\N_0$ be given by
\begin{equation*}
    a_{n}'(x)=\begin{cases}
        a_{n}(x) \quad \text{\rm if } \ a_{n}(x) \in \clP \, ,\\
        0 \qquad \quad \, \text{otherwise,}
    \end{cases}
\end{equation*}
where $\clP$ denotes the set of prime numbers. 
For any positive function $\vphi:\N\to [0, \infty)$, define the set
\[
\clE_{\ell}'(\vphi)
:=
\left\{ x\in [0,1): a_n'(x)\cdots a_{n+\ell-1}'(x)\geq \vphi(n) \text{ for i. m. } n\in\N\right\}.
\]
Studying the set of continued fractions with restrictions involving the set of prime numbers is not new, see for example \cite{MR4900954, MR1487636, SchindlerZweimuller2023}. In particular, Schindler and Zweim\"uller gave a Borel-Bernstein Theorem for $\clE_{1}'(\vphi)$ \cite[Theorem 2.1]{SchindlerZweimuller2023}.
We extend upon their work and completely determine the size of $\clE_{\ell}'(\varphi)$. In particular we calculate $\leb(\clE_{\ell}'(\varphi))$, see Theorem~\ref{THM:LEBMEAS}, and $\hdim \clE_{\ell}'(\varphi)$, Theorem~\ref{TEO:MAIN:HDIM}, thus providing the prime analogue of \cite[Theorem 1.7]{HuangWuXu2020}.

\subsection{The prime and almost prime zeta functions} \label{almost prime zeta function}
The most famously known Dirichlet series, the Riemann zeta function $\zeta$, is defined for $s>1$ as
\begin{equation*}
    \zeta(s):=\sum_{n\in\N} \frac{1}{n^{s}}\, .
\end{equation*}
The study of the Riemann zeta function is a vast area of mathematics that we simply cannot give due discussion in this article. See, for example, \cite{EdwardsRiemannzeta, TitchmarshRiemannzeta} for introductions to the topic and references within. Other well-known forms of Dirichlet series is the prime zeta function and the almost prime zeta function. Let $\ell \in \N$ and $\Omega(\ell)$ denote the set of $\ell$-primes. That is, the set of positive integers with at most $\ell$ (not necessarily distinct) prime factors. For example, $\Omega(1)=\clP$, and $2^{2}, 6 \in \Omega(2)$. For fixed $\ell\in\N$ the almost prime zeta function is defined as
 \begin{equation*}
     P_{\ell}(s;M):=\sum_{k\in \Omega(\ell): \, k \geq M}\frac{1}{k^{s}}\, .
 \end{equation*}
 In the special case $\ell=1$, $P_{1}(s;M)$ is the prime zeta function. Various results on the almost prime zeta function have been proven; see for example, \cite{frobergPrimezeta}. In particular, the relations
\begin{equation*}
    P(s;1)=\sum_{k\in \N}\frac{\mu(k)}{k}\log \zeta(ks)\, , \quad \text{ and } \quad  \zeta(s)=1+\sum_{\ell\in\N} P_{\ell}(s;1)\, ,
\end{equation*}
where $\mu$ is the M\"obius function, are well-known properties of such functions \cite{Glaisher}. \par
In \cite{SchindlerZweimuller2023}, Schindler and Zweim\"uller give an asymptotic for the tail of the prime zeta function at $s=2$ (Proposition \ref{Pr:SchindlerZweimuller}).
In this article, we determine the asymptotic behaviour of the tail of the almost prime zeta function (Theorem \ref{THM:PZF}).
Although the case $s=2$ of this statement is crucial in one of our main results in continued fraction theory, it is interesting on its own right.\par 
\subsection{Main results}{\label{S:MainResults}}
Before we state our results, we briefly summarise some frequently used notation. Throughout, $\clP$ is the set of prime numbers.  For any $\bfa=\sanu\in\N^{\N}$, and $r,s\in \N$ with $r\leq s$, write 
$\bfa[r,s]:=a_r a_{r+1} \ldots a_s$.
We extend the notation to finite words in an obvious way.
Next, for any $(x_n)_{n\geq 1}$ and $(y_n)_{n\geq 1}$ of positive real number, we write
 $x_n\ll y_n$ if for some constant $\kappa>0$ and every large $n$ we have $x_n\leq \kappa y_n$; if $\kappa$ depends on some parameter $\lambda$, we write $x_n\ll_{\lambda} y_n$.
Finally, we write $x_n\sim y_n$ when $\displaystyle\lim_{n\to\infty} \frac{x_n}{y_n}=1$.

Our first result describes the asymptotic behaviour of the Prime Zeta Function
 \begin{teo01}\label{THM:PZF}
For any $\ell\in\N$, real $s>1$, and every sufficiently large $M$ (depending on $\ell$) we have
\[
P_{\ell}(s;M)\asymp_{\ell}
 \frac{(\log\log M)^{\ell-1}}{M^{s-1}\log M}.
\]
\end{teo01}
We then use the case $s=2$ of Theorem \ref{THM:PZF} to calculate $\leb(\clE_{\ell}'(\vphi))$.
\begin{teo01}\label{THM:LEBMEAS}
For any $\vphi:\N\to[3, \infty)$, we have
\begin{equation*}
    \leb(\clE_{\ell}' (\vphi)) = \begin{cases}
        0 \quad \text{\rm if} \quad \displaystyle\sum_{n\geq 1} \frac{(\log\log \varphi(n))^{\ell-1}}{\varphi(n)\log\varphi(n)} < \infty \, , \\
        1 \quad \text{\rm if} \quad \displaystyle\sum_{n\geq 1} \frac{(\log\log \varphi(n))^{\ell-1}}{\varphi(n)\log\varphi(n)} = \infty \, .
    \end{cases}
\end{equation*}
\end{teo01}
Note that for $\ell=1$ we recover the Borel-Bernstein-type result by Schindler and Zweim{\"u}ller \cite[Theorem 2.1]{SchindlerZweimuller2023}. 
Also, for $\ell \geq 2$,  we cannot remove the restriction $\vphi(n)\geq 3$, $n\in\N$. 
Certainly, for the constant function $\vphi(n)=e$, $n\in\N$, the series in Theorem \ref{THM:LEBMEAS} is convergent. 
However, by virtue of Birkhoff's ergodic theorem, every finite block of natural numbers appears infinitely many times in the continued fraction expansion of almost every irrational $x\in (0,1)$, so $\leb(\clE_{\ell}'(\vphi))=1$, 
\par
On certain occasions, one would prefer to have the function $\varphi$ dependent on the size of the denominator of the convergent, $q_{n}$, rather than $n$. Consider the set 
\[
\clF_{\ell}(\vphi)
:=
\left\{ x\in [0,1): a_{n}'(x)\cdots a_{n+\ell-1}'(x)\geq \vphi(q_n)\text{ for i.m.} \  n\in\N\right\}.
\]
The exponential growth of $(q_n)_{n\geq 1}$ and the existence of the Lévy-Khinchin constant allow us to derive the next corollary.
The omitted proof follows almost verbatim that of \cite[Corollary 1.6]{HuangWuXu2020}.
\begin{coro01}\label{CORO:LEBMEAS}
For any $\vphi:\N\to [3, \infty)$,
\begin{equation*}
    \leb( \clF_{\ell} (\vphi) )=\begin{cases}
        0 \quad \text{\rm if} \quad \displaystyle\sum_{n\geq 1} \frac{(\log\log \varphi(n))^{\ell-1}}{n\varphi(n)\log\varphi(n)}<\infty\, ,\\
        1 \quad \text{\rm if} \quad \displaystyle\sum_{n\geq 1} \frac{(\log\log \varphi(n))^{\ell-1}}{n\varphi(n)\log\varphi(n)}=\infty\, .
    \end{cases}
\end{equation*}
\end{coro01}
Our last result tells us that the Hausdorff dimension of $\clE_{\ell}'(\vphi)$ is not affected by the primality restriction. 
\begin{teo01}\label{TEO:MAIN:HDIM}
For any function $\vphi:\N\to \mathbb{R}_{>0}$ and $\ell\in\N$, we have
\begin{equation}\label{Eq:HuangWuXuPrime}
\hdim \clE_{\ell}'(\vphi)
=
\hdim \clE_{\ell}(\vphi).
\end{equation}
\end{teo01}
An intermediate step in the proof of Theorem \ref{TEO:MAIN:HDIM} is a prime version of {\L}uczak's Theorem (Theorem \ref{TEO:LuczakTheorem}). 
For any function $\vphi:\Na\to\RE$, define
\[
\clE_{\ell}''(\vphi)
:=
\left\{x\in[0,1]: a_{n}'(x)\ldots a_{n+\ell-1}'(x)\geq \varphi(n) \, \text{ for all } \,  n\in\N\right\}. 
\]
\begin{teo01} \label{TEO:luczak for primes}
For any $b,c>1$, we have
\begin{equation*}
    \hdim \clE_{\ell}''(\vphi_{c,b})=\hdim \clE_{\ell}'(\vphi_{c,b})=\frac{1}{b+1}.
\end{equation*}
\end{teo01}

Our techniques suffice to show the prime analogue of a general theorem by Hussain and Shulga \cite{MR4912795}. 
To state it, for $\ell\in\N$ and $A_0,\ldots, A_{\ell-1}>1$, let $S_{\ell}'(A_0,\ldots, A_{\ell-1})$ be the set of irrational numbers $x\in [0,1]$ such that for infinitely many $n\in\N$ we have
\[
A_0^n\leq a_n'(x)\leq 2A_0^n, \;
A_{1}^n\leq a_{n+1}(x)\leq 2A_{1}^n, \;
\ldots,\;
A_{\ell-1}^n\leq a_{n+\ell-1}(x)\leq 2A_{\ell-1}^n.
\]
Also, let $S_{\ell}'(A_0,\ldots, A_{\ell-1})$ be the set of irrational set of irrational numbers $x\in [0,1]$ such that for infinitely many $n\in\N$ we have
\[
A_0^n\leq a_n'(x)\leq 2A_0^n, \;
A_{1}^n\leq a_{n+1}'(x)\leq 2A_{1}^n, \;
\ldots,\;
A_{\ell-1}^n\leq a_{n+\ell-1}'(x)\leq 2A_{\ell-1}^n.
\]

\begin{propo01}
For any $\ell\in\N$ and any $A_0,\ldots, A_{\ell-1}>1$, we have 
\[
\hdim S_{\ell}(A_0,\ldots, A_{\ell-1})
=
\hdim S_{\ell}'(A_0,\ldots, A_{\ell-1}).
\]
\end{propo01}

\noindent {\bf Acknowledgments.} This research is supported by the ARC Discovery Project 200100994. For part of this work Benjamin Ward was a Leverhulme Early Career Research Fellow.
We thank Cid Reyes Bustos, Tanja Schindler, Nikita Shulga, and Tim Trudgian for helpful discussions.

\subsection*{Structure of the paper.} 
In Section~\ref{Sec:Prelims}, we recall some preliminary results concerning the distribution of primes, Hausdorff dimension, and pressure functions and dimensional numbers.
In Section~\ref{Section:Proof:THM:PZF}, we show Theorem \ref{THM:PZF}.
Next, we prove Theorem~\ref{THM:LEBMEAS} in Section~\ref{Proof:THM:LEBMEAS}.
Section \ref{SECTION:TEO:luczak for primes} contains the proof of Theorem \ref{TEO:luczak for primes}.
In Section~\ref{Section:Proof:TEO:MAIN:HDIM}, we prove Theorem~\ref{TEO:MAIN:HDIM}.
Lastly, Section~\ref{Section:FurtherResearch} contains problems and suggestions for further research.

%
%
\section{Preliminaries}\label{Sec:Prelims}
\subsection{The distribution of primes}\label{primetheory} 
In this section, we recall two fundamental results on the distribution of primes: Mertens' Second Theorem and the Prime Number Theorem (PNT).
For a proof of the former, see, for example \cite[Theorem 2.7.(d)]{MontgomeryVaughan2007}.
A proof of the latter can be found in \cite[Section 6.2]{MontgomeryVaughan2007}.

\begin{teo01}[Second Mertens' Theorem] \label{second mertens theorem}
There is a constant $L$ such that
\[
\sum_{\substack{p\leq x  \\ p \in \clP}} \frac{1}{p} 
=
\log \log (x) + L + O\left( \frac{1}{\log x} \right)
\text{ as }
x\to\infty.
\]
\end{teo01}
For any $x>0$, write
\[
\pi(x)
: =
\#\left\{ p\in \clP: p\leq x\right\}.
\]
\begin{teo01}[Prime Number Theorem]\label{TEO:PNT}
We have $\pi(n)\sim \frac{n}{\log n}$.
\end{teo01}
The PNT allows us to estimate the number of primes in large intervals.
\begin{propo01}\label{CORO:PNT}
For any $\gamma>1$ and $n\in\Na$, let $c_n(\gamma)$ be determined by  
\[
\#\left(\clP\cap \left[ \gamma^n, 2\gamma^n\right]\right)
=
c_n(\gamma)
\frac{n\log \gamma}{\gamma^n}.
\]
Then, $c_n(\gamma)\to 1$ as $n\to \infty$.
\end{propo01}
Regarding primes in small intervals, we state a recent theorem by Bordignon, Johnston, and Starichkova.
\begin{propo01}{\cite[Lemma 15]{arxiv220709452}}\label{PR:ShortInterval}
If $x>e^{20}$, then $[0.999x,x)\cap \clP\neq \vac$.
\end{propo01}
A quantitative refinement of PNT was shown by Rosser in 1941.
\begin{propo01}[{\cite[Theorem 29]{Rosser1941}}] \label{Lem:Rosser}
If $x\geq 55$, then 
\[
\frac{x}{2+\log x}
<
\pi(x)
<
\frac{x}{-4+\log x}.
\]
\end{propo01}

The following proposition, which provides an asymptotic on the summation of the reciprocals of squared primes bounded from below, can be deduced from the PNT and Karamata's Theorem for slowly varying functions \cite[Theorem 0.6]{Resnick1987}. It was already noted by Schindler and Zweim\"{u}ller \cite[p. 267]{SchindlerZweimuller2023}. 
\begin{propo01} \label{Pr:SchindlerZweimuller}
For $M\in \N$ sufficiently large, we have that
\begin{equation*}
  \sum_{\substack{p>M \\ \, p \in \clP}}\frac{1}{p^{2}} \asymp \frac{1}{M\log M}\, .
\end{equation*}
\end{propo01}

\subsection{Continued fractions}\label{Subsection:contfractheory}
In this section, we recall classical results on continued fractions. 
Their proofs can be found in any standard reference on the topic, such as \cite{KhinchinCFBook} or \cite{HardyWright1979}.\par 
For any real number $y$, let $[y]$ be the integer part of $y$.
The Gauss map $T:[0,1)\to [0,1)$, which is intimately connected to continued fractions, is defined by
\begin{equation*}
\forall x\in [0,1), 
\quad
    T(x)
    =
    \begin{cases}
    x^{-1}-[x^{-1}], &\text{\rm if } \quad x\neq 0,\\
    0, &\text{\rm if }\quad  x= 0.
    \end{cases}
\end{equation*}
Let $a_1:(0,1)\to\N$ be given by $a_1(x)=[x^{-1}]$ and, for each $n\in\N_{\geq 2}$ and $x\in (0,1)$ such that $T^{n-1}(x)\neq 0$, put $a_n(x)=a_1(T^{n-1}(x))$. 
Then, we have 
\[
x
=
[a_1(x),a_2(x),a_3(x),\ldots] 
:=
\cfrac{1}{a_1(x) + \cfrac{1}{a_2(x) + \ddots}}.
\]
(see, for example, \cite[Chapter 2]{EinsiedlerWardBook}). 
For $n\in\N$, the \textit{$n$th convergent} of $x$ is the quotient
\[
\frac{p_n}{q_n}
:=
[a_1,\ldots, a_n],
\]
where $p_n, q_n\in\N$ are coprime. 
The sequences $(p_n)_{n\geq 1}$ and $(q_n)_{n\geq 1}$ satisfy 
\begin{align*}
p_{-1}=1, \, p_{0}=0,  \, &\text{ and } \, p_{n}=a_{n}p_{n-1}+p_{n-2} \,\text{ for }\,n\in\N,\\
q_{-1}=0, \,q_{0}=1,  \, &\text{ and } \, q_{n}=a_{n}q_{n-1}+q_{n-2} \,\text{ for } \,n\in\N.
\end{align*}
We can define $p_n$ and $q_n$ for a specific $n$-tuple of integers $\bfa$ with the same formula as above. 
In this case, we denote them by $p_n(\bfa)$ and $q_n(\bfa)$.\par

For $m, n\in\N$, $\bfa=(a_1,\ldots, a_n)\in\N^n$, and $\bfb=(b_1,\ldots,b_m)\in\N^n$, we write
\[
\bfa\bfb
:=
(a_1,\ldots, a_n,b_1,\ldots, b_m).
\]
\begin{propo01}\label{Pr:CF:EST}
Let $n\in \Na$ and $\bfa=(a_1,\ldots,a_n)\in\N^n$ be arbitrary.
\begin{enumerate}[i.]
\item \label{Pr:CF:EST:i} \cite[Lemma 2.1]{Wu2006} For any $k\in\{1,\ldots,n\}$, we have
\[
\frac{a_{k} + 1}{2}
\leq 
\frac{q_n(\bfa)}{q_{n-1}(a_1,\ldots,a_{k-1},a_{k+1},\ldots, a_n)}
\leq 
a_{k}+1.
\]
\item \label{Pr:CF:EST:ii} \cite[Proposition 1.1]{HensleyBook}
For any $m \in\N$ and $\bfb\in\N^m$, 
\[
q_{n}(\bfa)q_{m}(\bfb)
\leq
q_{n+m}(\bfa \bfb)
\leq 
2q_{n}(\bfa)q_{m}(\bfb).
\]
\end{enumerate}
\end{propo01}

For each $n\in\N$ and $\bfa=(a_1,\ldots, a_n)\in\N^n$, the \textit{fundamental interval} $I_n(\bfa)$ is defined as
\[
I_n(\bfa)
:=
\left\{ x\in [0,1): a_1(x)=a_1,\ldots, a_n(x)=a_n\right\}.
\]
\begin{propo01}[{\cite[Section 12]{KhinchinCFBook} }]\label{Pr:CF:FundamentalIntervals}
Let $n\in\N$ and $\bfa\in\N^n$ be arbitrary. 
\begin{enumerate}[i.]
\item If $p_n=p_n(\bfa)$, $q_n=q_n(\bfa)$, $p_{n-1}=p_{n-1}(a_1,\ldots, a_{n-1})$, $q_{n-1}=q_{n-1}(a_1,\ldots, a_{n-1})$, then
\[
I_n(\bfa)
=
\begin{cases}
    \left[ \displaystyle\frac{p_n}{q_n}, \frac{p_n+p_{n-1}}{q_n+q_{n-1}}\right), &\text{\rm if } n \text{ is even}, \\
    \left( \displaystyle\frac{p_n+p_{n-1}}{q_n+q_{n-1}},\frac{p_n}{q_n} \right], &\text{\rm if } n \text{ is odd}. 
\end{cases}
\]
\item There are some absolute constants for which the next estimate holds:
\[
\leb(I_n(\bfa))
\asymp
\frac{1}{q_n^2(\bfa)}\asymp_{n}\frac{1}{a_{1}^{2}\dots a_{n}^{2}}.
\]
\end{enumerate}
\end{propo01}

Proposition~\ref{Pr:CF:FundamentalIntervals} and the disjointedness of individual fundamental intervals imply the following result.

\begin{propo01}\label{PR:UnionFundamentalIntervals}
For $n\in\N$, $\bfa \in \N^{n}$, and $a,b\in\N$, we have
\begin{equation*}
\leb\left(
\bigcup_{a\leq j \leq b} I_{n+1}(\bfa,j)
\right)
=
    \frac{(b+1)-a}{(aq_n +q_{n-1})((b+1)q_n +q_{n-1})}\, .
\end{equation*}
\end{propo01}

\subsection{Hausdorff dimension}
We will provide a definition for the Hausdorff dimension of a set and some of its basic properties. 
For more details on this subject, we refer the reader to \cite{FalconerBook2014}. 
We begin by defining the Hausdorff Measure, recalling that a $\delta$-cover of a set $F\subseteq \RE$ is a countable collection of sets with diameter less than $\delta$ that cover $F$. 

\begin{def01}[{\cite[p. 44-45]{FalconerBook2014}}] 
    Suppose $F$ is a subset of $\RE^n$ and $ s \geq 0$. For each $\delta > 0$, we define 
    \begin{equation*}
        \clH_\delta^s(F) = \inf \left \{ \sum_{i=1}^\infty \diam(U_i)^s : \{U_i\} \  \text{is a}  \ \delta-\text{cover of} \ F  \right \} . 
     \end{equation*}
     Then the Hausdorff measure of the set is given by 
     \begin{equation}
         \clH^s(F) = \lim_{\delta \to 0} \clH^s_\delta(F). 
     \end{equation}
     Furthermore, the Hausdorff Dimension of a set $F$ is given by 
\begin{equation*}
    \hdim F = \inf \{ s \geq 0 : \clH^s(F) = 0\} = \sup\{ s: \clH^s(F) = \infty\}\, .
\end{equation*} 
\end{def01}

The following properties follow immediately from this definition. 
\begin{propo01}[{\cite[p.48]{FalconerBook2014}}]\label{Pr:HD:monotonicity} \label{Pr:HD:CountableStability}
If $ E \subset F$, then $\hdim E \leq \hdim F$. Furthermore, if $F_1, F_2, \ldots$ is a countable collection of subsets of $\RE$, then $$\hdim \bigcup_{i=1}^\infty F_i = \sup_{1 \leq i < \infty} \hdim F_i.$$ 
\end{propo01}

The following two results give us tools which can be used to estimate the lower bound of the Hausdorff dimension of certain sets. The first result is the Mass Distribution Principle (MDP). 
\begin{teo01}\label{THM:MDP}
Let $\mu$ be a mass distribution on a set $F\subset\R$ and suppose that for some $s>0$, there are numbers $c>0$ and $\varepsilon >0$ such that 
\begin{equation*}
    \mu(U) \leq c\diam(U)^s
\end{equation*}
for all intervals $U$ with $\diam(U)\leq \varepsilon$. Then $\mathcal{H}^s(F) \geq \mu(F)/c$ and $s \leq \hdim F$ . 

\end{teo01}

Lemma \ref{LE:FalconerExample4.6} below gives us a consequence of the Mass Distribution Principle when applied to certain Cantor sets (for a proof, see \cite[Example 4.6]{FalconerBook2014}).

\begin{lem01}\label{LE:FalconerExample4.6}
Let $(E_n)_{n\geq 1}$ be a sequence of finite unions of disjoint compact intervals contained in $[0,1]$ such that $E_n\supseteq E_{n+1}$ for each $n\in\N$. Suppose that each set $E_n$ is formed by at least $m_n\geq 2$ intervals with maximum widths of the intervals tending to zero as $n\to \infty$. Suppose the intervals for each $n$th level are separated by gaps of length at least $\veps_n>0$ with $(\veps_n)_{n\geq 1}$ decreasing. Then
\[
\hdim \bigcap_{n\in\N} E_n
\geq 
\liminf_{n\to\infty} \frac{\log(m_1\cdots m_n)}{-\log m_{n+1}\veps_{n+1}}.
\] 
\end{lem01}

\subsection{Pressure functions and dimensional numbers}\label{Subsec:PresureFunctions}
We use pressure functions to succinctly express the Hausdorff dimension of certain sets determined by continued fractions.
In this section, we introduce some notation and recall basic properties of pressure functions. 
We only consider pressure functions in the context of continued fractions and the Gauss map; for an overview of the general theory, we refer the reader to \cite{MR2003772}. 

Let $\vphi:[0,1] \to\RE$ be any function. 
For $n\in\N$ and $x\in [0,1]$, denote by $S_n(\vphi;x)$ the ergodic sum
\[
S_n(\vphi;x)
:=
\vphi(x) + \vphi(T(x)) + \ldots + \vphi\left(T^{n-1}(x)\right).
\]
For any non-empty subset $\clA\subseteq \Na$, write
\[
X_{\clA}
:=
\left\{ x\in [0,1): a_n(x)\in \clA \text{ for all } n\in\Na\right\}.
\]
It is well known  \cite[Lemma 2.1.2]{MR2003772} that the limit
\begin{equation}\label{Eq:Def:PressureFunction}
P_{\clA}(T,\vphi)
:=
\lim_{n\to\infty} 
\frac{1}{n} 
\log\left( 
\sum_{\bfa\in \clA^n} \sup_{x\in I_n(\bfa)\cap X_{ \clA}} \exp\left( S_n(\vphi,x)\right)
\right)
\end{equation}
exists. We call $P_{\clA}(T,\vphi)$ the \textit{pressure function} of the \textit{potential} $\vphi$. Define 
\[
\Var_n(\vphi)
:=
\sup\left\{ |\vphi(x)-\vphi(y)|:I_n(x)=I_n(y)\right\}.
\]
to be the \textit{$n$th variation} of $\vphi$. When evaluating the pressure function of the potential $\vphi$ the $n$th variation is particularly useful, as shown by the following proposition.

\begin{propo01}[{\cite[Proposition 2.4]{MR3162824}}]\label{Pr:PF:01}
Let $\vphi:[0,1]\to\RE$ be such that 
\[
\Var_1(\vphi)<\infty 
\;\text{ and }\; 
\lim_{n\to\infty} \Var_n(\vphi)=0.
\]
Then, the following statements are true:
\begin{enumerate}[1.]
    \item For every non-empty $\clA\subseteq \N$, the supremum in \eqref{Eq:Def:PressureFunction} can be replaced with the evaluation at any $x\in I_n(\bfa)\cap X_{\clA}$; 
    \item We have
    \[
P_{\Na}(T,\vphi)
=
\sup \left\{ P_{\clA}(T,\vphi): \clA\subseteq \Na\;\text{\rm is finite}\right\}.
\]
\end{enumerate}
\end{propo01}
Given any real $B>1$ and $f:[0,1]\to\RE$ continuous, define for each $s\in [0,1]$ the potential
\[
\vphi( f,s ;x)
:= 
-f(s)\log B - s\log|T'(x)|.
\]
Since for every $x\in [0,1)$ where $T^n$ is differentiable satisfies
\[
(T^n)'(x)
=
\frac{(-1)^n}{(xq_{n-1}-p_{n-1})^2},
\]
Proposition~\ref{Pr:PF:01} yields
\[
P_{\Na}(T,\vphi(f,s))
=
\lim_{n\to\infty} 
\frac{1}{n} 
\log\left( 
\sum_{\bfa\in A^n} \frac{1}{B^{nf(s)} q_{n}^{2s}(\bfa)} 
\right).
\]
We apply these observations on the sequence of continuous and non-decreasing functions $(f_{\ell})_{\ell \geq 1}$ defined in Theorem \ref{TEO:HWX:HDIM}.
For any $\ell,n\in\N$ and $\vac\neq \clA\subseteq \N$, define
\begin{align*}
    t_B^{(\ell)}(\clA)
    &:=
    \inf\{ s\geq 0: P_{\clA}(T;-f_{\ell } (s)\log B- s\log |T'|)\leq 0\}, \\
    t_B^{(\ell)}(\clA,n)
    &:=
    \inf\left\{ s\geq 0: \sum_{\bfa\in\clA^n} \frac{1}{B^{nf_{\ell }(s)}q_n^{2s}} \leq 0\right\}.
\end{align*}
Whenever $\clA=\{1,\ldots, M\}$ for some $M\in\N$, we write
\[
t_B^{(\ell)}(M):=t_B^{(\ell)}(\{1,\ldots, M\})
\quad\text{ and }\quad
t_B^{(\ell)}(M,n):=t_B^{(\ell)}(\{1,\ldots, M\},n).
\]
See \cite[Corollary 2.9 and Proposition 2.10]{HuangWuXu2020} for a proof of the next proposition.
\begin{propo01}\label{LE:DIMENSIONALNUMBER}
Let $\ell\in\N$ be arbitrary.
\begin{enumerate}[1.]
\item For every real $B>1$ we have
\[
t_B^{(\ell)}(\N)
=
\sup\left\{ t_B^{(\ell)}(\clA):\clA\subseteq \N \text{ is finite} \right\} 
=
\lim_{n\to\infty} t_B^{(\ell)}(\N,n) 
=
\lim_{M\to\infty} t_B^{(\ell)}(M).
\]
\item For every real $B>1$ and $M\in\Na$, 
\[
\lim_{n\to\infty} t_B^{(\ell)}(M,n)=t_B^{(\ell)}(M).
\]
\item  The next limits hold:
\[
\lim_{B\to 1} t_B^{(\ell)}(\N) =1 
\;\text{ and }\;
\lim_{B\to\infty } t_B^{(\ell)}(\N) = \frac{1}{2}.
\]
\end{enumerate}

\end{propo01}

\section{Proof of Theorem \ref{THM:PZF}}\label{Section:Proof:THM:PZF}
Let us recall part of Karamata's Theorem on regularly varying functions. For a proof and the full statement, see, for example, \cite[Theorem 0.6]{Resnick1987}.
\begin{propo01}\label{Prop:KaramataTheorem}
Let $F:\RE_{\geq 1} \to \RE_{\geq 0}$ be a measurable function for which there is some $\rho\in\RE$ such that
\[
\forall x>1
\quad
\lim_{t\to\infty}
\frac{F(tx)}{F(t)} = x^{\rho}.
\]
If $\rho<-1$, then $\int_{x}^{\infty} F(t)\md t<\infty$ and 
\[
\lim_{x\to\infty}
\frac{xF(x)}{\int_{x}^{\infty} F(t)\md t}= - \rho -1.
\]
\end{propo01}

\subsection{Proof of Theorem \ref{THM:PZF} }
For notation purposes for any positive integer $r\geq 1$ let
\begin{equation*}
    S(\ell,M,r,s):=\underset{ p_i \in \clP \, \, p_{i} \geq r}{\sum_{p_1 \dots p_\ell > M : \,}}\, \, \frac{1}{p_1^s\dots p_n^s}\, ,\quad \text{ and } \quad  
    \widehat{S}(\ell,M,r,s):=\underset{p_i \in \clP \, \, p_1\geq p_2 \geq \dots \geq p_\ell\geq r}{\sum_{p_1 \dots p_\ell > M : \,}}\frac{1}{p_1^s\dots p_\ell^s}\, .
\end{equation*}
Notice that for any $r\in \N$ 
\begin{equation*}
    \widehat{S}(\ell,M,r,s) \leq S(\ell,M,r,s) \leq \ell ! \widehat{S}(\ell,M,r,s) \, ,
\end{equation*}
 with the lower bound being trivial and the upper bound due to the fact that there are $\ell !$ ways to order $\ell$ prime numbers. Setting $r=1$, by our definition of $S(\ell,M,r,s)$, we will have completed the proof on showing that
 \begin{equation*}
     S(\ell,M,1,s)\asymp_{\ell} \frac{(\log\log M)^{\ell-1}}{M^{s-1}\log M}\, .
 \end{equation*}\par 
The proof is by induction. We begin by showing the base $\ell=1$ case.\par

Consider an arbitrary $s>2$ and define $F_s:[1,\infty)\to \RE$ by
\[
\forall x\in[1, \infty)
\quad
F_s(x) = \frac{1}{(x\log x)^s}.
\]
Then, for every $x>1$ we have
\[
\lim_{t\to\infty}
\frac{F_s(tx)}{F_{s}(t)} = x^{-s}.
\]
As a consequence, since $-s<-1$, Proposition \ref{Prop:KaramataTheorem} implies
\[
\lim_{x\to\infty} 
\frac{xF_s(x)}{\int_{x}^{\infty} F_s(t)\md t} = s- 1.
\]
Therefore, for large $x>1$ we have
\[
\int_{x}^{\infty} \frac{\md t}{(t\log t)^s}
\asymp_s
\frac{1}{x^{s-1}(\log x)^s}.
\]
We may now apply the PNT to conclude that for large $N$ we have
\begin{equation} \label{base case}
\sum_{n\geq N} \frac{1}{p_n^s}
\sim 
\sum_{n\geq N} \frac{1}{(n\log n)^s}
\asymp_s 
\frac{1}{N^{s-1}(\log N)^s}.
\end{equation}
Here and in the following equation only, the subscript on $p$ means that $p_{n}$ is the $n$th prime number. Later on $p_{n}$ will be used to denote our $n$th choice of prime number. Using \eqref{base case}, for $M$ sufficiently large, we have that

\begin{align*}
    S(1,M,1,s) 
    &=
    \sum_{p\in \clP : \, p> M}\frac{1}{p^{s}} 
    \asymp 
    \sum_{n>\pi(M)}\frac{1}{p_{n}^{s}} \\
    &\asymp_{s} \frac{1}{\pi(M)^{s-1}(\log \pi(M))^{s}} \asymp \frac{1}{\left(\frac{M}{\log M}\right)^{s-1}(\log M)^{s}}=\frac{1}{M^{s-1}\log M}
\end{align*}

completing the base $\ell=1$ case. \par

Now, let us assume that the result is true for $\ell$ and show that it holds for $\ell +1$. We split the inductive step in two, giving an inductive argument for the upper bound and then an inductive argument for the lower bound. The observation that, for any $r\geq 1$
 \begin{equation} \label{prime count observation}
     S(\ell,M,r,s)=\sum_{p\in\clP: \, p\geq r}\frac{1}{p^{s}}S(\ell-1,\tfrac{M}{p},r,s) \quad \text{ and } \quad \widehat{S}(\ell,M,r,s)=\sum_{p\in\clP: \, p\geq r}\frac{1}{p^{s}}\widehat{S}(\ell-1,\tfrac{M}{p},p,s)
 \end{equation}
 is useful in both cases. Our assumption on $M$ being sufficiently large is to guarantee that we can apply Theorem~\ref{second mertens theorem} and Proposition~\ref{Pr:SchindlerZweimuller}. We also require in the $\ell$-induction assumption that $M^{\tfrac{1}{\ell+1}}$ is sufficiently large so that
\[
\sum_{\substack{p_1\cdots p_{\ell}\geq M^{\tfrac{1}{1+\ell}}\\ p_1,\ldots,p_{\ell}\in\clP}} \frac{1}{p_1^s\cdots p_{\ell}^s} 
\asymp_{\ell}
\frac{(\log\log M)^{\ell-1}}{M^{(s-1)\tfrac{1}{\ell+1}}\log M}.
\] 
 Beginning with the upper bound, for $\ell+1$, we have that 
\begin{equation*}
    S(\ell+1,M,1,s)\asymp_{\ell} \widehat{S}(\ell+1,M,1,s)=\sum_{p_{\ell+1} \in \clP}\frac{1}{p_{\ell+1}^{s}} \widehat{S}\left(\ell,\tfrac{M}{p_{\ell+1}},p_{\ell+1},s\right)\, . 
\end{equation*}
We split this into two summations depending on the size of $p_{\ell+1}$. Namely,
\begin{equation*}
    \widehat{S}(\ell+1,M,1,s)= \sum_{p_{\ell+1}\leq M^{\tfrac{1}{\ell+1}}: \, p_{\ell+1} \in \clP}\frac{1}{p_{\ell+1}^{s}} \widehat{S}\left(\ell,\tfrac{M}{p_{\ell+1}},p_{\ell+1},s\right)\, + \sum_{p_{\ell+1}>M^{\tfrac{1}{\ell+1}}: \, p_{\ell+1} \in \clP}\frac{1}{p_{\ell+1}^{s}} \widehat{S}\left(\ell,\tfrac{M}{p_{\ell+1}},p_{\ell+1},s\right)\, .
\end{equation*}
Begin with the first sum. Notice that $p_{\ell+1}\leq M^{\tfrac{1}{\ell+1}}$ implies that 
\begin{equation} \label{bound1}
    \frac{M}{2} \geq \frac{M}{p_{\ell+1}}\geq M^{1-\tfrac{1}{\ell+1}} \geq M^{\tfrac{1}{\ell+1}} \quad \text{ and so } \log \tfrac{M}{p_{\ell+1}}\asymp_{\ell} \log M \, .
\end{equation}
So
\begin{align}
    \sum_{p_{\ell+1}\leq M^{\tfrac{1}{\ell+1}}: \, p_{\ell+1} \in \clP}\frac{1}{p_{\ell+1}^{s}} \widehat{S}\left(\ell,\tfrac{M}{p_{\ell+1}},p_{\ell+1},s\right) &\asymp_{\ell} \sum_{p_{\ell+1}\leq M^{\tfrac{1}{\ell+1}}: \, p_{\ell+1} \in \clP}\frac{1}{p_{\ell+1}^{s}} S\left(\ell,\tfrac{M}{p_{\ell+1}},p_{\ell+1},s\right) \label{upper bound step 1}\\
    &\leq \sum_{p_{\ell+1}\leq M^{\tfrac{1}{\ell+1}}: \, p_{\ell+1} \in \clP}\frac{1}{p_{\ell+1}^{s}} S\left(\ell,\tfrac{M}{p_{\ell+1}},1,s\right) \label{less than condition}\\
    &\asymp \sum_{p_{\ell+1}\leq M^{\tfrac{1}{\ell+1}}: \, p_{\ell+1} \in \clP}\frac{1}{p_{\ell+1}^{s}} \frac{\left(\log \log \tfrac{M}{p_{\ell+1}}\right)^{\ell-1}}{\left(\tfrac{M}{p_{\ell+1}}\right)^{s-1}\log \tfrac{M}{p_{\ell+1}}}, \nonumber
\end{align}
where the last step follows from our inductive assumption, which is applicable by our lower bound in \eqref{bound1}. Notice that it is only step \eqref{less than condition} that prevents us from obtaining an asymptotic statement. Using that $\log\tfrac{M}{p_{\ell+1}}\asymp_{\ell}\log M$ we have that
\begin{equation} \label{upper bound step 2}
    \sum_{p_{\ell+1}\leq M^{\tfrac{1}{\ell+1}}: \, p_{\ell+1} \in \clP}\frac{1}{p_{\ell+1}^{s}} \widehat{S}\left(\ell,\tfrac{M}{p_{\ell+1}},p_{\ell+1}\right) \ll_{\ell} \sum_{p_{\ell+1}\leq M^{\tfrac{1}{\ell+1}}: \, p_{\ell+1} \in \clP}\frac{1}{p_{\ell+1}} \frac{\left(\log \log M\right)^{\ell-1}}{M^{s-1}\log M} \, , 
\end{equation}
and so, by applying Second Mertens Theorem (Lemma~\ref{second mertens theorem}) we have that
\begin{equation} \label{sum 1 outcome}
    \sum_{p_{\ell+1}\leq M^{\tfrac{1}{\ell+1}}: \, p_{\ell+1} \in \clP}\frac{1}{p_{\ell+1}^{s}} \widehat{S}\left(\ell,\tfrac{M}{p_{\ell+1}},p_{\ell+1}\right) \ll_{\ell} \frac{\left(\log \log M\right)^{\ell}}{M^{s-1}\log M} \, .
\end{equation}
For the second part of the summation, since we have that $p_{\ell+1}>M^{\tfrac{1}{\ell+1}}$ and the fact that our primes our bounded from below by $p_{\ell+1}$, we can write
\begin{equation*}
    p_{1},p_{2},\dots,p_{\ell+1}>M^{\tfrac{1}{\ell+1}} \, .
\end{equation*}
So, multiplying together, we have that $p_{1}\ldots p_{\ell+1}>M$. Thus the restriction on the product of the primes becomes redundant, leading us to
\begin{align}
    \sum_{p_{\ell+1}>M^{\tfrac{1}{\ell+1}}: \, p_{\ell+1} \in \clP}\frac{1}{p_{\ell+1}^{s}} \widehat{S}\left(\ell,\tfrac{M}{p_{\ell+1}},p_{\ell+1},s\right)&\asymp_{\ell} \sum_{p_{\ell+1}>M^{\tfrac{1}{\ell+1}}: \, p_{\ell+1} \in \clP}\frac{1}{p_{\ell+1}^{s}} S\left(\ell,\tfrac{M}{p_{\ell+1}},p_{\ell+1},s\right) \nonumber\\
    & = \left(\sum_{p>M^{\tfrac{1}{\ell+1}}: \, p \in \clP}\frac{1}{p^{s}} \right)^{\ell+1} \nonumber\\
    &\overset{\text{case $\ell=1$}}{\asymp_{\ell}} \left(\frac{1}{M^{(s-1)\tfrac{1}{\ell+1}}\log M^{\tfrac{1}{\ell+1}}}\right)^{\ell+1} \nonumber\\
    &\asymp_{\ell} \frac{1}{M^{s-1}(\log M)^{\ell+1}}\, .\label{sum outcome 2}
\end{align}
Combining \eqref{sum 1 outcome} and \eqref{sum outcome 2} we get our required upper bound. That is,
\begin{equation*}
    S(\ell+1,M,1,s)\ll_{\ell} \frac{\left(\log \log M\right)^{\ell}}{M^{s-1}\log M} + \frac{1}{M^{s-1}(\log M)^{\ell+1}} \ll  \frac{\left(\log \log M\right)^{\ell}}{M^{s-1}\log M} \, .
\end{equation*}
In order to show the lower bound for $\ell+1$, using the previous observation \eqref{prime count observation} we have that
\begin{equation} \label{lower bound split}
S(\ell+1,M,1,s)
\geq 
\sum_{\substack{p\leq M^{\tfrac{1}{\ell+1}}\\ p\in\clP}} 
\frac{1}{p^s}S(\ell,\tfrac{M}{p},1,s) \, .
\end{equation}
This is similar to our calculation in the upper bound, except this time we no longer have step \eqref{less than condition} since our prime numbers are already bounded from below by $1$. So following steps \eqref{upper bound step 1},\eqref{upper bound step 2} and \eqref{sum 1 outcome} we get the asymptotic
\begin{equation*}
    \sum_{\substack{p\leq M^{\tfrac{1}{\ell+1}}\\ p\in\clP}} 
\frac{1}{p^s}S(\ell,\tfrac{M}{p},1,s) \asymp_{\ell} \frac{(\log\log M)^{\ell}}{M^{s-1}\log M}
\end{equation*}
Hence combining with \eqref{lower bound split} we obtain 
\begin{equation*}
    S(\ell+1,M,1,s) \gg_{\ell} \frac{(\log\log M)^{\ell}}{M^{s-1}\log M} 
\end{equation*}
completing the proof.

\section{Lebesgue Measure. Proof of Theorem \ref{THM:LEBMEAS}}{\label{Proof:THM:LEBMEAS}}
In \cite{KleinbockWadleigh2018}, Kleinbock and Wadleigh showed Lemma \ref{Le:KW} below using mixing the properties of the Gauss map and Philipp's quantitative Borel-Cantelli Lemma.
Their result is stated in terms of the Gauss measure $\mu$, $\md \mu = \frac{1}{(\log 2)(1+x)} \md \leb$, but we can obviously write it in terms of the Lebesgue measure.

\begin{lem01}[{\cite[Lemma 3.5]{KleinbockWadleigh2018}}] \label{Le:KW}
Let $r\in \N$ and $(B_{n})_{n\geq 1}$ be a sequence of countable unions of fundamental intervals of level $r$. If
\begin{equation*}
A:=\left\{ x \in [0,1) : T^{n}(x)\in B_{n} \ \text{\rm for i.m. } n \in \N\right\} \, ,
\end{equation*}
then
\begin{equation*}
    \leb(A)=\begin{cases}
        0, &\text{\rm if } \quad \displaystyle\sum_{n=1}^{\infty}\leb(B_{n})< \infty \, , \\
        1, &\text{\rm if } \quad \displaystyle\sum_{n=1}^{\infty}\leb(B_{n})= \infty \, . 
    \end{cases}
\end{equation*}
 \end{lem01}
 
Let $T\in\N$ be such that the conclusion of Theorem~\ref{THM:PZF} holds for all $M\geq T$.
We further assume
\begin{equation} \label{inital assumption}
    \text{there exists $n_{0}\in\N$ such that  $\varphi(n)>M$ for all $n\geq n_{0}$,}
\end{equation}
and then show how this condition can be removed at the end of the proof.
\par

The first step is to rewrite $\clE^{'}_{\ell}(\varphi)$ in the form of the set $A$ in Lemma~\ref{Le:KW}. To do this, note that 
\begin{align*}
    \clE^{'}_{\ell}(\varphi)&=\limsup_{n\to \infty} \left\{ x \in [0,1): \, a_{n}'(x)\ldots a_{n+\ell-1}'(x)\geq \varphi(n) \right\}\, ,\\
   &= \limsup_{n\to \infty} \left\{ x \in [0,1): \, a_{n+1}'(x)\ldots a_{n+\ell}'(x)\geq \varphi(n+1) \right\}\,
\end{align*}
Then, in terms of the $n$th iterate of the Gauss map, and using that $a_{n}(x)=a_{1}\left(T^{n-1}(x)\right)$, we have
\begin{align*}
    \clE^{'}_{\ell}(\varphi)&=\limsup_{n\to \infty} \left\{ x \in [0,1): \, a_{1}'\left(T^{n}(x)\right)\ldots a_{\ell}'\left(T^{n}(x)\right)\geq \varphi(n+1) \right\}\, ,\\
    &=\limsup_{n\to \infty} \left\{ x \in [0,1): \, T^{n}(x) \in \left\{y\in[0,1): a_{1}'\left(y\right)\ldots a_{\ell}'\left(y\right)\geq \varphi(n+1)\right\} \right\} \, .
\end{align*}
For each $n\in\N$, put
\begin{equation*}
    B_{n}:=\left\{y\in[0,1): a_{1}'\left(y\right)\ldots a_{\ell}'\left(y\right)\geq \varphi(n+1)\right\}\, .
\end{equation*}
 The sets $B_{n}$ can be written as a countable collection of disjoint level $\ell$ fundamental intervals:
\begin{equation*}
    B_{n}=\underset{p_{1}\ldots p_{\ell}\geq \varphi(n+1)}{\bigcup_{p_{1},\ldots,p_{\ell} \in \clP}} I_{\ell}(p_{1},\ldots,p_{\ell} )\, .
\end{equation*}
The second step is to calculate the Lebesgue measure of each $B_{n}$.
By Proposition~\ref{Pr:CF:FundamentalIntervals}, 
\begin{equation*}
   \leb\left( I_{\ell} (p_{1},\ldots,p_{\ell}) \right)\asymp_{\ell}\frac{1}{p_{1}^{2}\ldots p_{\ell}^{2}}\, ,
\end{equation*}
so
\begin{equation*}
    \leb(B_{n})=\underset{p_{1}\ldots p_{\ell}\geq \varphi(n+1)}{\sum_{p_{1},\ldots,p_{\ell} \in \clP }} \frac{1}{p_{1}^{2}\ldots p_{\ell}^{2}}\, .
\end{equation*}
Since $\varphi(n)\geq T$ for $n\geq n_0$, Theorem~\ref{THM:PZF} yields
\begin{equation*}
    \leb(B_{n}) \asymp_{\ell} \frac{(\log\log \varphi(n+1) )^{\ell-1}}{\varphi(n+1)\log\varphi(n+1)}\, . 
\end{equation*}
Therefore, in view of Lemma~\ref{Le:KW}, 
\begin{equation*}
    \leb(\clE_{\ell}' (\vphi)) = \begin{cases}
        0 \quad \text{\rm if} \quad \displaystyle\sum_{n\geq 1} \frac{(\log\log \varphi(n))^{\ell-1}}{\varphi(n)\log\varphi(n)} < \infty \, , \\
        1 \quad \text{\rm if} \quad \displaystyle\sum_{n\geq 1} \frac{(\log\log \varphi(n))^{\ell-1}}{\varphi(n)\log\varphi(n)} = \infty \, ,
    \end{cases}
\end{equation*}
completing the proof under our initial assumption. \par

Now, to show the theorem without assumption \eqref{inital assumption}, suppose there is a strictly increasing sequence $(n_{i})_{i\geq 1}$ such that $\varphi(n_{i})\leq T$ for all $i \in \Na$; in particular,
\begin{equation}
\forall i\in\N
\quad
    \frac{\left(\log\log\varphi(n_{i})\right)^{\ell-1}}{\varphi(n_{i})\log\varphi(n_{i})}
    \geq 
    F(T):=
    \frac{\left(\log\log T\right)^{\ell-1}}{T\log T}\, .
\end{equation}
We are, then, in the divergence case:
\begin{equation*}
    \sum_{n\geq 1} \frac{(\log\log \varphi(n))^{\ell-1}}{\varphi(n)\log\varphi(n)} 
    > 
    \sum_{i\geq 1} \frac{(\log\log \varphi(n_{i}))^{\ell-1}}{\varphi(n_{i})\log\varphi(n_{i})} 
    \geq 
    \sum_{i\geq 1} F(T) = \infty \, .
\end{equation*}
Define $\psi:\Na \to \Na_{\geq 2}$ by
\begin{equation*}
\forall n\in\N\quad
    \psi(n)=\max\{T,\varphi(n)\}\, .
\end{equation*}
Then $\clE^{'}_{\ell}(\psi)\subseteq \clE^{'}_{\ell}(\varphi)$. 
Moreover, $\psi(n)\geq T$ for all $n\in\N$, so
\begin{equation} \label{psi divergence}
    \leb(\clE^{'}_{\ell}(\psi))=1 \quad \text{\rm if } \quad \sum_{n\geq 1} \frac{(\log\log \psi(n))^{\ell-1}}{\psi(n)\log\psi(n)}=\infty \, .
\end{equation}
To complete the proof, notice that 
\begin{equation*}
    \sum_{n\geq 1} \frac{(\log\log \psi(n))^{\ell-1}}{\psi(n)\log\psi(n)} 
    >
    \sum_{i\geq 1} \frac{(\log\log \psi(n_{i}))^{\ell-1}}{\psi(n_{i})\log\psi(n_{i})}
    =
    \sum_{i\ge 1} F(T)=\infty\,.
\end{equation*}
Therefore, by \eqref{psi divergence} and $\clE^{'}_{\ell}(\psi)\subseteq \clE^{'}_{\ell}(\varphi)$, we conclude
$
    \leb(\clE^{'}_{\ell}(\varphi)) \geq \leb(\clE^{'}_{\ell}(\psi))=1.
$
Since we were in the divergence case, the proof is complete.
\qed
\section{Hausdorff dimension. Proof of Theorem \ref{TEO:luczak for primes}}{\label{SECTION:TEO:luczak for primes}}
In this section, we prove Theorem \ref{TEO:luczak for primes}. 
Fix $c,b>1$ and recall that $\varphi_{b,c}(n)=c^{b^{n}}$ for $n\in\N$. 
For any $\ell\in\N$, the upper bound for $\hdim \clE_{\ell}'(\vphi_{b,c})$ is immediate from Theorem \ref{TEO:LuczakTheorem} and $\clE_{\ell}'(\vphi_{b,c})\subseteq \clE_{\ell}(\vphi_{b,c})$. 
For the lower bound, we use Lemma~\ref{LE:FalconerExample4.6} (cfr.\cite{MR1464376}).
Note that we need to ensure our intervals are non-empty; that is, we need to check that the primes are dense enough `fill out' the levels. 

\subsection{Case $\ell=1$.}
As mentioned, we only need to show:
\begin{equation}\label{Eq:LuLowerELL1}
    \hdim \clE_{1}''(\vphi_{b,c}) \geq \frac{1}{b+1}.
\end{equation}
We briefly recall the notation from Lemma~\ref{LE:FalconerExample4.6}:
\begin{enumerate}[i.]
    \item $E_{k}$ is the union of $k$ level basic interval of the fractal set $\bigcap_{k \in \Na} E_{k}$.
    \item $m_{k}$ is a lower bound on the children of any $(k-1)$th level interval.
    \item $\varepsilon_{k}$ is a lower bound in the gap between any two $k$th level intervals.
\end{enumerate}
First, for each $k\in\N$, consider
\begin{equation*}
    \clI_{k}:=\prod_{j=1}^{k} \left(\clP \cap \left[c^{b^{j}},3c^{b^{j}}\right]\right)\, ,
\end{equation*}
which is non-empty by Bertrand's Postulate.
Define the set
\begin{equation*}
    E_{k}:=\bigcup_{\mathbf{d}=(d_{1},\ldots,d_{k})\in \clI_{k}} I_{k}(\mathbf{d})\, ,
\end{equation*}
hence
\begin{equation*}
    \bigcap_{k\in\N} E_k
    \subseteq 
    \clE_{1}''(\varphi_{c,b})\subseteq \clE_{1}'(\varphi_{c,b})\, .
\end{equation*}\par
Next, we look at the integers $m_{k}\geq 2$. 
Note that for $k\in\N$ and $\bfd\in \clI_k$, we have $I_{k}(\mathbf{d})\subset I_{k-1}(\mathbf{f})$ if and only if $\mathbf{d}= \mathbf{f} P$ for some prime $P$ with $ c^{b^{k}}\leq P \leq 3c^{b^{k}}$, so we need a lower bound on the number of primes in $[c^{b^{k}},3c^{b^{k}} ]$. 
For all large $k$, Lemma~\ref{Lem:Rosser} gives
\begin{align*}
\#\left(\clP\cap\left[c^{b^{k}},3c^{b^{k}}\right]\right) &=\pi\left(3c^{b^{k}}\right) - \pi\left(c^{b^{k}}\right)\\
&> \frac{3c^{b^{k}}}{2+ \log 3c^{b^{k}}}
-
\frac{c^{b^{k}}}{-4+ \log c^{b^{k}}} \\
&>
\frac{3c^{b^{k}}}{2\log c^{b^{k}}} 
- 
\frac{c^{b^{k}}}{\log c^{b^{k}}} \\
&=
\frac{c^{b^{k}}}{2b^{k}\log c } =: m_{k}.
\end{align*}
We now define $\veps_k$. Take $\bfd\in\clI_{k-1}$. 

For any two primes $c^{b^{k}}<P<R<3c^{b^{k}}$ and $k$ sufficiently large we trivially have $R-P\geq 2$, so between the fundamental intervals $I_{k}(\bfd \,P)$ and $I_k(\bfd\,R)$ is at least one fundamental interval, which has length at least
\begin{align*}
\leb\left( I_{k}\left(\bfd \, \left(\left\lfloor3c^{b^k}\right\rfloor+1\right) \right) \right) &=\left(q_{k}\left(\bfd \, \left(\left\lfloor3c^{b^k}\right\rfloor+1\right) \right) \left(q_{k}\left(\bfd \, \left(\left\lfloor3c^{b^k}\right\rfloor+1\right) \right)+q_{k-1}(\bfd)\right)\right)^{-1}\\
&\geq 
\left( 2 q_{k}^2\left(\bfd \, \left(\left\lfloor3c^{b^k}\right\rfloor+1\right) \right) \right)^{-1} \\
&\geq \left( 2\prod_{i=1}^{k-1}(d_{i}+1)^{2}\left(\left\lfloor3c^{b^k}\right\rfloor+1\right)^{2} \right)^{-1}
\end{align*}
by Proposition~\ref{Pr:CF:FundamentalIntervals} and a repeated application of Proposition~\ref{Pr:CF:EST}.\ref{Pr:CF:EST:i}.

Since
\begin{equation*}
    2\prod_{i=1}^{k-1}(d_{i}+1)^{2}\left(\left\lfloor3c^{b^k}\right\rfloor+1\right)^{2} \leq 2^{2k+3}9^{k}c^{2b+2b^{2}+\ldots + 2b^{k-1}+2b^{k}}\leq 36^{k+1} c^{2\frac{b^{k+1}-b}{b-1}},
\end{equation*}
we may pick
\[
\veps_k
:=
36^{-k-1} c^{-2 \frac{b^{k+1}-b}{b-1}}.
\]
Then, as $k\to\infty$, we have
\begin{align*}
\frac{ \log(m_1\cdots m_{k-1})}{-\log (m_k\veps_k)}
&=
\frac{(b + b^2 + \ldots + b^{k-1})\log c + O(k^2)}{-b^k\log c + 2 \frac{b^{k+1}-b}{b-1} \log c + O(k)} \\
&=
\frac{b^k \frac{\log c}{b-1} (1+o(1))+O(k^2) }{ \frac{2b^{k+1} - 2b - b^{k+1} + b^k}{b-1} \log c + O(k)} \\
&=
\frac{ 1+o(1) }{ b+1 +  o(1)} \sim \frac{1}{b+1}.
\end{align*}
We now conclude \eqref{Eq:LuLowerELL1} from Lemma~\ref{LE:FalconerExample4.6}.

\subsection{General case.}
Take $\ell\in\N_{\geq 2}$.
For the lower bound, note that $\clE_{1}''(\varphi_{b,c}) \subseteq \clE_{\ell}''(\varphi_{b,c})\subseteq \clE_{\ell}^{'}(\varphi_{b,c})$, so
\begin{equation} \label{upper bound}
    \hdim \clE_{\ell}'(\varphi_{b,c})
    \geq 
    \hdim \clE_{\ell}''(\varphi_{b,c}) 
    \geq 
    \hdim \clE_{1}''(\varphi_{b,c})=\frac{1}{b+1}\, .
\end{equation}
As mentioned earlier, the upper bound is obvious from Theorem \ref{TEO:HWX:HDIM}.
Nevertheless, we provide a proof without relying on this result.
Indeed, take any $y \in \clE_{\ell}^{'}(\varphi_{b,c})$. Then there exists an infinite sequence $(n_{j})_{j\geq 1}$ such that $a_{n_{j}}^{'}(y)\ldots a_{n_{j}+\ell-1}^{'}(y)\geq c^{b^{n_{j}}}$ for all $j\in \N$. By the geometric mean, we have that
\begin{equation} \label{maximum term}
    \max_{0\leq t\leq \ell-1} a_{n_{j}+t}^{'}(y) 
    \geq c^{\tfrac{b^{n_{j}}}{\ell}}\, .
\end{equation}
Let us define for $t\in \{0,1,\ldots, \ell-1\}$ the numbers
\begin{equation*}
    c_{t}:=c^{\tfrac{1}{\ell b^{t}}}>1\, , \quad\text{ and so }\, \, c_{t}^{b^{n+t}}=c^{\tfrac{b^{n}}{\ell}}\, .
\end{equation*}
Since the sequence $(n_{j})_{j\in\N}$ is infinite, there exists $0\leq v\leq \ell-1$ such that \eqref{maximum term} occurs at $a_{n_{j_{i}}+v}'(y)$ for an infinite sequence $(n_{j_{i}})_{i\geq 1}$. Hence, for all $i\in\N$,
\begin{equation*}
    a_{n_{j_{i}}+v}'(y)\geq c^{\frac{b^{n_{j_{i}}}}{\ell}}=c_{v}^{b^{n_{j_{i}}+v}}\, ,
\end{equation*}
so $y\in \clE_{1}^{'}(\varphi_{b,c_{v}})$. This shows that
\begin{equation*}
    \clE_{\ell}'(\varphi_{b,c})
    \subseteq \bigcup_{0\leq v\leq \ell-1}\clE_{1}' (\varphi_{b,c_{v}})\,,
\end{equation*}
Then, by the countable stability of the Hausdorff dimension (Proposition \ref{Pr:HD:CountableStability}),
\begin{equation} \label{lower bound}
    \hdim \clE_{\ell}^{'}(\varphi_{b,c}) \leq \hdim \bigcup_{0\leq v\leq \ell-1}\clE_{1}^{'}(\varphi_{b,c_{v}}) =\max_{0\leq v\leq \ell-1}\hdim \clE_{1}^{'}(\varphi_{b,c_{v}})=\frac{1}{b+1}\, .
\end{equation}
Combining \eqref{upper bound} and \eqref{lower bound}, we have that
\begin{equation*}
    \frac{1}{b+1}
    \leq 
    \hdim \clE_{\ell}''(\varphi_{b,c}) 
    \leq 
    \hdim \clE_{\ell}'(\varphi_{b,c}) \leq \frac{1}{b+1}\, .
\end{equation*}
This completes the proof for general $\ell$.
\qed
\section{Hausdorff dimension. Proof of Theorem \ref{TEO:MAIN:HDIM}}\label{Section:Proof:TEO:MAIN:HDIM}
Fix a natural number $\ell$ and a function $\vphi:\N\to \RE_{>0}$. 
By $\clE_{\ell}'(\vphi)\subseteq \clE_{\ell }(\vphi)$, the result holds if $\leb(\clE_{\ell}(\vphi))>0$. 
If $\leb(\clE_{\ell}'(\vphi))=0$ and $\leb(\clE_{\ell}(\vphi))=1$, then
\[
\sum_{n\geq 1} \frac{(\log\log \varphi(n))^{\ell-1}}{\varphi(n)\log\varphi(n)} < \infty 
\;\text{ and }\;
\sum_{n=1}^{\infty}\frac{\left(\log \varphi(n)\right)^{\ell-1}}{\varphi(n)}=\infty
\]
by Theorem \ref{THM:LEBMEAS} and Theorem \ref{huang wu xu lebesgue theory} respectively.
From the divergence of the second series, we may show that $B_{\vphi}=1$ so, by the discussion to come, $\hdim \clE_{\ell}'(\vphi)=1$, proves the result.
Thus, we may assume that both $\clE_{\ell}'(\vphi)$ and $\clE_{\ell}(\vphi)$ are Lebesgue-null sets.
By the monotonicity of the Hausdorff dimension, Theorem \ref{TEO:MAIN:HDIM} will be proven once we show
\begin{equation}\label{EQ:HWX:PF01}
\hdim \clE_{\ell }'(\vphi)
\geq 
\hdim \clE_{\ell}(\vphi).
\end{equation}

As usual, most of the work occurs when $1<B_{\vphi}<\infty$. 
To solve this case, we first consider a restricted family of functions.

\subsection{A particular case}
Fix a real number $B>1$.
In this section, we apply the Mass Distribution Principle on a Cantor set $E_B\subseteq  \clE_{\ell}'(\vphi_B)$ to show
\begin{equation}\label{eq:HD of EB}
\hdim \clE_{\ell}'(\vphi_B)
\geq 
\hdim \clE_{\ell}(\vphi_B).
\end{equation}

\subsubsection{Setting up the parameters}
We now establish the parameters needed in the overall construction:
\begin{enumerate}[i.]
\item Let $s$ and $\delta$ be positive real numbers such that 
\[
\frac{1}{2}<s<t_B^{(\ell)}
\;\text{ and }\;
\frac{1}{2}<s-2\delta < s.
\]
\item Define the real numbers $\alpha_0,\ldots, \alpha_{\ell - 2}>1$ by
\begin{equation}\label{EQ:DefAlphas}
\forall j\in\{0,\ldots, \ell-2\}
\quad
\alpha_j=B^{\frac{s^{\ell-1-j}(2s-1)(1-s)^j}{s^{\ell }- (1-s)^{\ell}}}.
\end{equation}

\item Let $M,N\in\N$ be so large that the next conditions hold:
\[
s<t_B(M,N;\ell); 
\]
\begin{equation}\label{EQ:PF:THMHD:N}
N> 
\max\left\{
e^{20},
\frac{5}{s\delta} +1, 
\frac{2 \ell }{\delta},
\frac{2 \ell \log 2}{\delta\log \alpha_0}
\right\};
\end{equation}
for each $i\in\{0,\ldots, \ell-2\}$, if $c_n(\alpha_i)$ is as in Corollary \ref{CORO:PNT}, we have
\[
\forall n\in\Na_{\geq N}
\quad
c_n(\alpha_i)<2;
\]
and
\begin{equation}\label{Eq:BndDeltaN}
\frac{N^{\ell} 2^{\ell} \log \alpha_0\cdots \log \alpha_i}{(\alpha_0\cdots \alpha_i)^{\delta N} }
=
\frac{2^{\ell}N^{\ell} \log\alpha_0 \cdots \log \alpha_{i-2}}{(\alpha_0(\alpha_1\cdots \alpha_{i })^2\alpha_{i+1})^{s\delta N} }
< 
1.
\end{equation}
The equality follows from Lemma \ref{Le:alphaIdentities} below.

\item Let $(l_j)_{j\geq 1}$ be a rapidly increasing sequence of natural numbers, say $ l_{n+1}> 2 l_n$ for all $n\in\N$. 
\item Let $(n_j)_{j\geq 0}$ be the sequence in $\N$ given by
\[
n_0 = -(\ell-1) 
\quad\text{ and }\quad
\forall j\in\N \quad n_{j}: = n_{j-1} + \ell +l_{j}N .
\]
\end{enumerate}
\begin{lem01}[{\cite[Lemma 5.1]{HuangWuXu2020}}]\label{Le:alphaIdentities}
The next estimates hold:
\begin{align*}
\forall j\in\{0,1,\ldots, \ell - 3 \}
\quad
\frac{1}{\alpha_0 \cdots \alpha_{j} }
&=
\left(
\frac{1}{\alpha_0 (\alpha_1\cdots \alpha_{j})^2 \alpha_{j+1} } 
\right)^s, \\
\frac{1}{\alpha_0 \cdots \alpha_{\ell-2}}
&=
\left(
\frac{1}{B \alpha_1\cdots \alpha_{\ell-2} } 
\right)^s, \\
B \alpha_0^s 
&\geq B^{2s}. 
\end{align*}
\end{lem01}
A word of warning, we consider blocks of length $\ell$ while the proof in \cite{HuangWuXu2020} considers blocks of length $k+1$.
\subsubsection{A Cantor Set}
If $\ell\geq 2$, let $E_B$ be the set of irrational numbers $x\in [0,1)$ such that
\begin{enumerate}[i.]
\item For every $j\in\N$, we have 
\[\forall i\in\{0,1,\ldots,\ell-2\}
\quad
\alpha_{i}^{n_j}\leq a_{n_j+i}'(x)\leq 2 \alpha_{i}^{n_j}.
\]
\item For every $j\in\N$, we have
\[
\left(\frac{B}{\alpha_0\cdots \alpha_{\ell -2}}\right)^{n_j}
\leq 
a_{n_j+\ell-1}'(x)
\leq
2\left(\frac{B}{\alpha_0\cdots \alpha_{\ell - 2}}\right)^{n_j}.
\]
\item For every other $n\in\N$, we have $1\leq a_n(x)\leq M$.
\end{enumerate}
When $\ell=1$, we define $E_B$ as the set of irrational numbers $x\in [0,1]$ with the next properties:
\begin{enumerate}[1.]
\item For all $j\in\N$, we have $B^{n_j}\leq a_{n_j}'(x) \leq 2 B^{n_j}$, and
\item For all $n\in \N\setminus \{n_j:j\in\N\}$, we have $a_n(x)\leq M$.
\end{enumerate}
\subsubsection{Construction of the measure}
Assume $\ell\geq 2$. For $j\in\N$ and $i\in\{0,\ldots, \ell-2\}$, define
\[
\clP_{i,j}
:=
\clP\cap [\alpha_i^{n_j}, 2 \alpha_i^{n_j}]
\quad\text{ and }\quad
c_{i,j}:=c_{n_j}(\alpha_i)
\]
(see Proposition \ref{CORO:PNT} for how $c_{n_{j}}(\alpha_{i})$ is defined); when $i=\ell-1$, we define
\[
\clP_{\ell-1,j}
:=
\clP
\cap 
\left[ \left(\frac{B}{\alpha_0\cdots \alpha_{\ell - 2}}\right)^{n_j}, \,2 \left(\frac{B}{\alpha_0\cdots \alpha_{\ell - 2}}\right)^{n_j}\right]
\quad\text{ and }\quad
c_{\ell-1, j}:=c_{n_j}\left(\left(\frac{B}{\alpha_0\cdots \alpha_{\ell - 2}}\right)^{n_j}\right).
\]
When $\ell=1$, for $j\in\N$ we define 
\[
\clP_{j}
:=
\clP\cap [B^{n_j}, 2 B^{n_j}].
\]
Let $D$ be the set of sequences of partial quotients of elements in $E_B$. That is, 
\[
D:=
\left\{
\sanu\in\N^{\N} : [a_1,a_2,a_3,\ldots]\in E_B
\right\} \, .
\]
For $n\in\N$, define
\[
D_n 
:=
\left\{
\bfa[1,n]: \bfa\in D
\right\}
\]
For $n\in\N$ and $\bfa\in D_n$, the \emph{fundamental set} $J_n(\bfa)$ is
\[
J_n (\bfa)
=
\bigcup_{\substack{a_{n+1}\in \N\\ \bfa a_{n+1}\in D_{n+1}}} \overline{I}_{n+1}(\bfa a_{n+1}).
\]
For any $\bfa=\sanu\in\N^{\N}$, $j\in\N$, and $i\in \{0,1,\ldots, l_j-1\}$, let $\bfb_j\in\Na^{l_jN}$ and $\bfb_{j,i}\in\N^{N}$ be determined by
\[
\bfb_{j}
=
\bfb_{j,1} \cdots \bfb_{j, l_j}
=
\bfa[ n_{j-1}+\ell,  n_j-1].
\]
This way, we have
\begin{align*}
\bfa
&=
\bfb_1 a_{n_1}\ldots a_{n_1 + \ell -1} \bfb_2 a_{n_2}\ldots a_{n_2 + \ell-1}\bfb_3 a_{n_3}\ldots a_{n_3 + \ell -1}\ldots\\
&=
\bfb_{1,1}\ldots \bfb_{1,l_1-1} a_{n_1}\ldots a_{n_1 + \ell-1} \bfb_{2,1}\ldots \bfb_{2,l_1-1} a_{n_2}\ldots a_{n_2 + \ell-1} \bfb_{3,1}\ldots \bfb_{3,l_3-1} a_{n_3}\ldots a_{n_3 + \ell-1}\ldots.
\end{align*}
We extend this notation for finite sequences in the obvious manner.

Let $u$ be given by
\[
u
:=
\sum_{\bfb\in \{1,\ldots,M\}^N} 
\left(\frac{1}{\alpha_0^Nq_N^2(\bfb)}\right)^s
=
\sum_{\bfb\in \{1,\ldots,M\}^N} 
\frac{1}{q_N^{2s}(\bfb)B^{Nf_{\ell}(s)} }
\geq 1.
\]
(for a proof of the second equality, see (5.2) in \cite{HuangWuXu2020}). \par
We define a measure $\mu$ on $E_B$ by defining it on the fundamental sets. 
First, consider $n\in\N$ such that $1\leq n < n_1 + \ell$ and $\bfa\in D_n$.
\begin{enumerate}[i.]
\item If $n=n_1-1=Nl_1$, put
\[
\mu\left( J_{n_1 - 1}(\bfa)\right)
=
\prod_{i=1}^{l_1} \frac{1}{u} \left( \frac{1}{ \alpha_0^N q_N^2 (\bfb_{1,j} )} \right)^s.
\]
\item If $n=n_1+i$ for some $i\in\{0,1,\ldots, \ell-1 \}$, put
\[
\mu\left( J_{n_1+ i}(\bfa)\right)
=
\frac{1}{\#\clP_{0,1} \cdots \#\clP_{i, 1} }
\mu\left( J_{n_1-1}(\bfb_1  )\right).
\]
\item If $n < n_1-1$, put
\[
\mu\left( J_n(\bfa)\right)
=
\sum_{\substack{\bfc \in \N^{n_1-1-n} \\ \bfa\bfc \in D_{n_1-1}}}
\mu\left(J_{n_1-1}(\bfa\bfc)\right).
\]
\end{enumerate}
Assume that we have defined $\mu$ on every fundamental set of order $n\leq n_{k} + \ell-1$ for some $k\in \N$. Take $n_{k} + \ell  \leq n \leq n_{k+1} + \ell -1$ and $\bfa\in D_n$.
\begin{enumerate}[i.]
\item  If $n=n_{k+1}-1$, define
\[
\mu\left( J_{n_{k+1}-1}(\bfa)\right)
=
\mu\left( J_{n_{k} +\ell }\left( \bfa[1,n_k+\ell -1]   \right)\right)
\prod_{j=1}^{l_{k+1}} 
\frac{u^{-1}}{ \alpha_0^{Ns} q_N^{2s}( \bfb_{k+1,j}  ) }.
\]
\item  If $n=n_{k+1}+i$ for some $i\in\{0,1\ldots, \ell-1 \}$, write
\[
\mu\left( J_{n_{k+1}+i}(\bfa)\right)
=
\frac{1}{\#\clP_{0, k+1}\,\cdots\, \#\clP_{i, k+1} }
\;
\mu\left( J_{n_{k+1}-1}(\bfa[1, n_{k+1}-1] )\right).
\]

\item If $n_{k}+\ell   \leq n \leq n_{k+1} - 2$, define
\[
\mu\left(J_n(\bfa)\right)
=
\sum_{\substack{\bfc\in \N^{n_{k+1}-1-n} \\ \bfa\bfc\in D_{n_{k+1}-1}}}
\mu\left(J_{n_{k+1}-1}(\bfa\bfc)\right)\, .
\]
\end{enumerate}
When $\ell=1$, we define $\mu$ as follows. 
Take any $n\in\Na$ and $\bfa\in D_n$.
\begin{enumerate}[i.]
\item If $n=n_1-1$, put
\[
\mu\left(J_{n_1-1}(\bfa\right))
=
\left( \frac{1}{B^{n_1-1}q_{n_1-1}^2(\bfa)}\right)^{ t_B^{(1)}(M,n_1-1) }.
\]
\item If $n=n_1$, put
\[
\mu\left(J_{n_1}(\bfa)\right)
=
\frac{1}{\# \clP_1}
\mu\left(J_{n_1-1}(\bfa[1,n_1-1]  \right).
\]
\item For $1\leq n \leq n_1-2$, put
\[
\mu\left(J_{n}(\bfa)\right)
=
\sum_{\bfc\in \{1,\ldots,M\}^{n_1-1-n} } \mu\left(J_{n_1-1}(\bfa \bfc)\right).
\]
\end{enumerate}

Assume that we have defined $\mu$ for all fundamental intervals of level $1,2,\ldots, n_{k}-1$ for some $k\in\Na$. Assume that $n_{k}\leq n\leq n_{k+1}-1$.
\begin{enumerate}[i.]
\item If $n=n_k$, then
\[
\mu\left(J_{n_k}(\bfa)\right)
=
\frac{1}{\#\clP_k} \mu\left( J_{n_k-1}( \bfa[1,n_k-1] \right).
\]
\item If $n= n_{k+1}-1$, then
\[
\mu\left(J_{n_{k+1}-1}(\bfa)\right)
=
\mu\left(J_{n_{k}}( \bfa[1,a_{n_k}] )\right)
\left(\frac{1}{B^{s_{n_{k+1} -n_k -1} } q^2_{m_{k+1}}(\bfb_{k+1}) } \right)^{t_B^{(1)}(M,n_{k+1} -n_k -1) }.
\]
\item If $n_{k}+1\leq n \leq n_{k+1} -2$, then
\[
\mu\left(J_{n}(\bfa)\right)
=
\sum_{\bfc\in \{1,\ldots, M\}^{n_{k+1}-n}} \mu\left(J_{n_{k+1}-1}(\bfa \bfc)\right).
\]
\end{enumerate}

In what follows, we assume $\ell\geq 2$ since the proof for $\ell=1$ follows closely that in \cite{WangWu2008}.
In both cases, the primality restriction gives rise to similar complications.

\subsection{Length and gap estimates}  
From Proposition \ref{PR:ShortInterval}, we know that for small $0<\lambda<10^{-4}$ and every large number $x$, the intervals $((1-\lambda)x,x)$ and $(x, (1-\lambda)^{-1}x)$ contain a prime. 
From this observation, the estimates in Subsection \ref{Subsection:contfractheory}, and standard computations, we conclude the next result (Proposition \ref{PR:ShortInterval} is used in items \ref{Le:LengthEst:ii} and \ref{Le:LengthEst:iii}, those involving the primality restriction).
\begin{lem01}\label{Le:EstDiamFundSets:HWX}
Let $n\in\N$ and $\bfa\in D_n$ be arbitrary and  $k\in\N_0$ such that $n_{k} + \ell -1 \leq n < n_{k+1} +\ell -1$.
The following estimates hold for some absolute constants.
\begin{enumerate}[i.]
\item \label{Le:LengthEst:i} If $n_{k} + \ell -1 \leq n < n_{k+1} -1$, then
\[
\diam( J_n(\bfa))
\asymp 
\frac{1}{q_n^2(\bfa)}.
\]

\item \label{Le:LengthEst:ii} If $n=n_{k+1}-1$, then 
\[
\diam(J_n(\bfa))
\asymp
\frac{ 1  }{  \alpha_0^{n_{k+1}} q_{n_{k+1} -1}^2(\bfa)}.
\]

\item \label{Le:LengthEst:iii} If $n_{k+1} \leq n \leq n_{k+1} +\ell -3$, then, for $j = n  - n_{k+1}$,
\[
\diam(J_n(\bfa))
\asymp
\frac{ 1 }{  \alpha_{j+1}^{n_{k+1}} ( \alpha_0 \cdots \alpha_j)^{ 2n_{k+1}} q_{n_{k+1}-1}^2(\bfa[1,n_{k+1}-1]) } .
\]
\item \label{Le:LengthEst:iv} If $n = n_{k+1} + \ell -2$, then 
\[
\diam(J_n(\bfa))
\asymp
\frac{ 1 }{ ( B\alpha_0\cdots \alpha_{\ell-2} )^{ n_{k+1}} q_{n_{k+1}-1}^2(\bfa[1,n_{k+1}-1]) }.
\]
\end{enumerate}

\end{lem01}

As an example of one of the calculations, suppose $n=n_{k+1}+j$ for $1\leq j<\ell-2$. The others follow similarly. By definition
\begin{align*}
    \diam(J_n(\bfa))
    &=
    \diam\left(\underset{a_{n+1} \in \clP}{\bigcup_{\alpha_{j}^{n_{k+1}}\leq a_{n+1}\leq 2\alpha_{j}^{n_{k+1}}}} I_{n+1}(\bfa a_{n+1})\right)\\
    &=
    \diam\left({\bigcup_{q_1\leq a_{n+1}\leq q_2}} I_{n+1}(\bfa a_{n+1})\right) \qquad \qquad \qquad\qquad \qquad\, \, \left(\text{$q_{1},q_{2}\in \clP$ solve }\begin{array}{c}
    q_{1}=\displaystyle\min_{p\in \clP}\{p\geq \alpha_{j}^{n_{k+1}}\} \\
    q_{2}=\displaystyle \max_{p\in \clP}\{p\leq 2\alpha_{j}^{n_{k+1}}\}
    \end{array} \right)\\
    &\asymp \frac{q_{2}-q_{1}}{\alpha_{j}^{2n_{k+1}}(\alpha_{0}\cdots\alpha_{j})^{2n_{k+1}}q_{n_{k+1}-1}^{2}(\bfa[1,n_{k+1}-1])} \qquad \text{(by Proposition~\ref{Pr:CF:EST}-\ref{Pr:CF:FundamentalIntervals})}\\
    &\asymp \frac{1}{\alpha_{j}^{n_{k+1}}(\alpha_{0}\cdots\alpha_{j})^{2n_{k+1}}q_{n_{k+1}-1}^{2}(\bfa[1,n_{k+1}-1])} \, \, \qquad \text{(by Proposition~\ref{PR:ShortInterval}).}
\end{align*}

For $n\in\N$ and $\bfa\in D_n$, let $G_n(\bfa)$ denote the minimal gap between $J_n(\bfa)$ and the nearest fundamental set of level $n$. 
In general, our fundamental sets do not coincide with the fundamental intervals defined in \cite{HuangWuXu2020}.
Our sets, however, are contained in the latter, so we can use the following gap estimate.

\begin{propo01}[{\cite[Lemma 5.3]{HuangWuXu2020}}]
For any $n\in\N$ and $\bfa\in D_n$, we have 
\[
G_n(\bfa) \geq \frac{1}{8M}\diam(I_n(\bfa)).
\]
\end{propo01}

\subsection{H\"older exponents for fundamental sets}
Now, we estimate the measure $\mu$ of the fundamental intervals in terms of their diameter.
\begin{propo01}\label{PROP:HoldExpFundSets}
For all $n\in\Na$ and $\bfa\in D_n$, we have
\[
\mu\left(J_n(\bfa)\right)
\ll_{\ell} 
 \diam\left(J_n(\bfa)\right)^{s(1-\delta)-\delta}.
\]
\end{propo01}
\begin{proof}
Take $n\in\N$ and $\bfa\in D_n$. 
The implied constants are either absolute or depend on $\ell$.
\begin{enumerate}[i.]
\item \label{case.i} Case $n=Nl$ for some $l\in\{0,\ldots, l_1-1\}$. From \eqref{EQ:PF:THMHD:N}, we get $q_{lN}(\bfa)^2 \geq 2^{lN-1} \geq 2^{\frac{4l}{s\delta}}$; hence,
\begin{align*}
\mu\left( J_n(\bfa) \right)
&\leq
\prod_{k=1}^{l} \left( \frac{1}{ \alpha_0^N q_N^2( \bfb_{1,k} ) } \right)^s \\
&<
\prod_{k=1}^{l} \left( \frac{1}{ q_N^2( \bfb_{1,k} ) }\right)^s \\
&\leq 
4^{l}
 \left( \frac{1}{ q_{lN}^2(\bfa  )}\right)^s \\
 &\leq 
 \left( \frac{1}{ q_{lN}^2(\bfa)}\right)^{s(1-\delta)} \\
 &\ll 
\diam\left(J_n(\bfa)\right)^{s(1-\delta)}.
\end{align*}
\item \label{case.ii} Case $n=n_1-1=l_1N$. 
\begin{align*}
\mu\left( J_n(\bfa)\right) 
&\leq 
\prod_{j=1}^{l_1} \left( \frac{1}{\alpha_0^N q_N^2(\bfb_{1,j})} \right)^s \\
&<
\left(\frac{1}{ \alpha_0^{n_1(1- 1/N)}  } \right)^s
\left( \frac{1}{q_n^2(\bfa)} \right)^{s(1-\delta)} \\
&\leq 
\left( \frac{1}{\alpha_0^{n_1}q_n^2(\bfa)} \right)^{s(1-\delta)}\\ 
&\ll
\diam\left(J_n(\bfa)\right)^{s(1-\delta)}.
\end{align*}

\item \label{case.iii} If $n\in \{n_1,\ldots, n_1 + \ell- 3 \}$, consider $i=n - n_1  \in \{0,\ldots, \ell -3\}$, then
\begin{align*}
\mu\left( J_n(\bfa)\right)
&=
\frac{1}{\#\clP_{0,1} \cdots \#\clP_{i, 1} }
\mu\left( J_{n_1-1}(\bfb_{1})\right) \\
&= 
\frac{c_{1,1}\cdots c_{i,1} n_1^{i} \log\alpha_0 \cdots \log \alpha_{i}}{(\alpha_0\cdots \alpha_{i})^{n_1} }
\mu\left( J_{n_1-1}(\bfb_{1})\right) \\
&\ll 
\frac{ 1 }{(\alpha_0\cdots \alpha_{i})^{n_1(1-\delta)} }
\left( \frac{1}{\alpha_0^{n_1}q_n^2(\bfb_1)} \right)^{s(1-\delta)}  &&\text{ (by \eqref{Eq:BndDeltaN})} \\
&\leq 
\frac{ 1 }{(\alpha_0(\alpha_1\cdots \alpha_{i})^2\alpha_{i+1} )^{sn_1(1-\delta)} }
\left( \frac{1}{\alpha_0^{N}q_n^2(\bfb_1)} \right)^{s(1-\delta)} && \text{(by Lemma \ref{Le:alphaIdentities})} \\
&\leq 
\left( \frac{1}{ \alpha_0^{N+n_1}(\alpha_1\cdots \alpha_{i})^{2n_1} \alpha_{i+1}^{n_1} q_n^2(\bfb_1)} \right)^{s(1-\delta)}  \\
&\ll
\diam(J_{n_1+i}(\bfa))^{s(1-\delta)}.
\end{align*}

\item \label{case.iv} If $n=n_1 + \ell - 2$, then
\begin{align*}
\mu\left( J_n(\bfa)\right)
&=
\mu(J_{n_1-1}(\bfb_1) )
\frac{1}{\#\clP_{0,1} \cdots \#\clP_{\ell-1, 1} } \\
&<
\left( \frac{1}{ \alpha_0^{n_1} q_{n_1-1}^2(\bfb_1) } \right)^{s(1-\delta)}  
\left( \frac{1}{\alpha_0 \cdots \alpha_{\ell-1}} \right)^{ n_1(1-\delta)} \\
&=
\left( \frac{1}{ \alpha_0^{n_1} q_{n_1-1}^2(\bfb_1) } \right)^{s(1-\delta)} 
\left( \frac{1}{ B\alpha_1 \cdots \alpha_{\ell-1} } \right)^{ n_1s(1-\delta)} \\
&=
\left( \frac{1}{ (B\alpha_0\alpha_1 \cdots \alpha_{\ell-1})^{n_1} q_{n_1-1}^2(\bfb_1) } \right)^{s(1-\delta)} \\
&\ll
\diam(J_{n_1+\ell-1}(\bfa))^{s(1-\delta)}\\
&=
\diam\left(J_n(\bfa)\right)^{s(1-\delta)}.
\end{align*}

\item \label{case.v} If $n=n_1+ \ell -1$, then
\begin{align*}
\mu\left(J_{n_1+\ell-1}(\bfa)\right)
&<
\frac{1}{B^{n_1(1-\delta)}}
\left(\frac{1}{\alpha_0^{n_1}q_{n_1-1}^2(\bfb_1)}\right)^{s(1 - \delta) } \\
&<
\left(\frac{1}{ B^{2n_1}  q_{n_1-1}^2(\bfb_1)}\right)^{s(1-\delta)} \\
& \ll 
\left(\frac{1}{ q_{n_1+\ell}^2(\bfa)}\right)^{s(1 - \delta)} \\
&\ll 
\diam\left(J_n(\bfa)\right)^{s(1 - \delta)} .
\end{align*}

\item \label{case.vi} If $N(l-1) < n < N l$ for some $l\in \{1,2,\ldots, l_1\}$, we use Case $a_{k+3}^1$ in \cite{HuangWuXu2020} with obvious modifications to conclude
\[
\mu\left( J_n(\bfa) \right)
\ll 
\diam(J_n(\bfa))^{s(1-\delta)}.
\]

\item \label{case.vii} If $n=n_{k+1}-1$ for some $k\in \N$, then
\begin{align*}
\mu\left( J_{n_{k+1} - 1}(\bfa)\right)
&=
\mu\left( J_{n_{k}+\ell}(\bfa[1,n_k+\ell] )\right)
\prod_{i=1}^{l_{k+1}} \frac{1}{u}\left( \frac{1}{\alpha_0^N q_N^2(\bfb_{k+1,i} )}\right)^s \\
&\leq 
\frac{1}{B^{n_{k+1}(1-\delta)}}
\mu\left( J_{n_{k} - 1}(\bfa[1,n_k-1])\right)
\prod_{i=1}^{l_{k+1}} \frac{1}{u}\left( \frac{1}{\alpha_0^N q_N^2(\bfb_{k+1,i} )}\right)^s \\
&\leq 
\left[ 
\prod_{t=1}^{k} 
\left(
\frac{1}{B^{n_t(1-\delta)}}
\prod_{i=1}^{l_t} \frac{1}{u}\left( \frac{1}{\alpha_0^N q_N^2(\bfb_{t,i} )}\right)^s
\right)
\right]
\prod_{i=1}^{l_{k+1}} \frac{1}{u}\left( \frac{1}{\alpha_0^N q_N^2(\bfb_{k+1,i} )}\right)^s
\end{align*}
For each $t\in\{1,\ldots, k\}$, we argue as before to get
\begin{align*}
\frac{1}{B^{n_t(1-\delta)}}
\prod_{i=1}^{l_t} \frac{1}{u}\left( \frac{1}{\alpha_0^N q_N^2(\bfb_{t,i} )}\right)^s
&\leq 
\frac{1}{B^{n_t(1-\delta)}}
\left( \frac{1}{ \alpha_0^{n_t} q_{Nl_t}^2(\bfb_t) } \right)^{s(1-\delta)} \\
&\leq 
\left( \frac{1}{B^{n_t} q_{Nl_t}^2(\bfb_t)}\right)^{s(1-\delta)} \\
&\leq 
4^{\ell+1}
\left(
\frac{1}{q_{l_tN+\ell+1}^2(\bfb_t,a_{n_t},\ldots, a_{n_t+\ell})}
\right)^{s(1-\delta)}. 
\end{align*}
We also have
\[
\prod_{i=1}^{l_{k+1} } \frac{1}{u}\left( \frac{1}{\alpha_0^N q_N^2(\bfb_{k+1,i} )}\right)^s
\leq 
\left( \frac{1}{\alpha_0^{n_j}q_{l_{k+1}N}^2(\bfb_{k+1})}\right)^{s(1-\delta)}.
\]
Also, recalling \eqref{EQ:PF:THMHD:N}, we have
\[
q_n^2(\bfa)
\geq
2^{n-1}
\geq 
4^{\frac{(k+1)N}{2} }
\geq 
4^{\frac{(k+1)(\ell+1)}{2} \frac{N}{\ell+1} }
\geq 
4^{\frac{(k+1)(\ell+1)}{\delta}  }
\]
These estimates yield
\begin{align*}
\mu\left(J_{n_{k+1}-1}(\bfa)\right)
&\leq
4^{\ell+1}4^{k+1}
\left( \frac{1}{\alpha_0^{n_{k+1} }q_{n_{k+1}-1}^2(\bfa)}\right)^{s(1-\delta)} \\
&\leq
4^{(k+1)(\ell+1)}
\left( \frac{1}{\alpha_0^{n_{k+1} }q_{n_{k+1}-1}^2(\bfa)}\right)^{s(1-\delta)} \\
&\leq
\left( \frac{1}{\alpha_0^{n_{k+1} }q_{n_{k+1}-1}^2(\bfa)}\right)^{s(1-\delta)-\delta}.
\end{align*}

\item For the remaining $n\in\N$, we simply combine Case \ref{case.vii} with the argument used in one of the previous six cases.
For example, when $n= n_{k+1} + \ell -1$, we use cases \ref{case.vii} and \ref{case.v}.

\end{enumerate}
\end{proof}

With the H{\"o}lder exponents of the fundamental sets and the gap estimates at hand, we may apply the argument in \cite[Subsection 5.4.4]{HuangWuXu2020} with virtually no changes to obtain the exponents for small balls.
\begin{propo01}
There is a positive constant depending on $M$ and $\ell$ such that for every sufficiently small $r>0$, for every $x\in E_B$ we have
\[
\mu(B(x,r))
\ll
r^{s(1-\delta)-\delta}.
\]
\end{propo01}

\subsection{Conclusion of the proof }
The above discussion and the Mass Distribution Principle imply that
\[
s(1-\delta)-\delta
\leq 
\hdim E_B
\leq 
\hdim \clE_{\ell}'(\vphi_B).
\]
We, thus, get \eqref{eq:HD of EB} by letting $\delta\to 0$ and $s\to t_B^{(\ell)}$.

\subsection{A slight extension}
In the above construction, we can specify where the blocks of large prime partial quotients may appear without affecting the dimension.
More precisely, for any $\vphi:\N\to\RE_{>0}$ and any infinite $\clN\subseteq\N$ define
\[
\clE_{\ell}'(\vphi;\clN)
=
\{x\in [0,1): a_{n}'(x)\cdots a_{n+\ell-1}'(x)\geq \vphi(n) \text{ for i.m. } n\in\clN\}.
\] 
\begin{lem01}\label{LE:HDIM:EPHIN}
If $1<B<\infty$ and $\clN\subseteq \N$ is infinite, then
\[
\hdim \clE_{\ell}'(\vphi_B;\clN)
=
\hdim \clE_{\ell}'(\vphi_B).
\]
\end{lem01}
\begin{proof}
We only discuss the lower bound, since $\clE_{\ell}'(\vphi_B,\clN) \subseteq \clE_{\ell}'(\vphi_B)$. 
Choose integers $M,N$ such that \eqref{EQ:PF:THMHD:N} holds.
We may further assume that $\clN=\{n_1<n_2<\ldots\}$ satisfies
\[
\forall j\in\N
\quad
\frac{n_{j+1}}{n_j} >1.
\]
Define $n_0=-\ell$ and let $(l_j)_{j\geq 1}$ and $(r_j)_{j\geq 1}$ be the sequences of non-negative integers given by
\[
\forall j\in\N
\quad
n_{j} - 1 - n_{j-1} - \ell    
 =    
  l_{j} N + r_{j} 
  \;\text{ and }\;
   0\leq r_{j}\leq N-1.  
\]
 Let $\alpha_0=\alpha_0(B)$, $\ldots$, $\alpha_{\ell-1}=\alpha_{\ell-1}(B)$ be as in \eqref{EQ:DefAlphas}.
Let $E_{\ell}(B;\clN)$ be the set of irrational numbers $x\in [0,1)$ with the next properties:
 \begin{enumerate}[1.]
\item For every $j\in\Na$, 
\begin{enumerate}[i.]
\item For all $i\in\{0,1,\ldots,\ell-2\}$,
\[
\alpha_i^{n_j} \leq a_{n_j+i}'(x)\leq 2 \alpha_i^{n_j};
\]
\item We have
\[
\left( \frac{B}{\alpha_0 \cdots \alpha_{\ell-2}}\right)^{n_j}  \leq a_{n_j+\ell-1}'(x) \leq  2\left( \frac{B}{\alpha_0 \cdots \alpha_{\ell-2}}\right)^{n_j}. 
\]
\item If $r_j>0$ and $n\in\N$ is such that $n_j + \ell + l_{j}N+1\leq n \leq n_{j+1}-1$, then $a_n(x)=2$.
\end{enumerate}

\item  All the other partial quotients $a_n(x)$ satisfy $a_n(x)\leq M$.
\end{enumerate}
   
Next, we define a measure $\mu$ just as in the proof of \eqref{EQ:HWX:PF01} except for the numbers $n\in\Na$ for which $n_j<n<n_{j+1}$, $r_j>0$ and $n_j + \ell + l_{j+1}N+1\leq n \leq n_{j+1}-1$ for some $j\in\N$. 
In such cases, for $\bfa\in D_n$,
\[
\mu\left( J_n(\bfa) \right)
=
\mu\left( J_{n_j + \ell+1 + l_{j+1}N}(a_1,\ldots, a_{n_j + \ell+1 + l_{j+1}N})\right).
\]
The argument we used to compute $\hdim E_B$ holds in our context too, hence $\hdim \clE_{\ell}(\vphi_B;\clN)=\hdim \clE_{\ell}(\vphi_B)$.
\end{proof}

\subsection{Proof of Theorem \ref{TEO:MAIN:HDIM}}
First, assume that  $B_{\vphi}=1$.
For any $C>1$, define the infinte set
\[
\clN_{C}
:=
\left\{ n\in\Na: \frac{\log \vphi(n)}{n} <\log C\right\}
\]
Then, $\clE_{\ell}'(\vphi_C;\clN_{C}) \subseteq \clE_{\ell}'(\vphi)$ and, by Lemma \ref{LE:HDIM:EPHIN}, 
\[
\hdim \clE_{\ell}(\vphi_C)
=
\hdim \clE_{\ell}'(\vphi_C;\clN_{C})
\leq 
\hdim \clE_{\ell}'(\vphi). 
\]
We conclude the result letting $C\to B_{\vphi}$ and using Proposition~\ref{LE:DIMENSIONALNUMBER}. 

Now assume that $B_{\vphi}=+\infty$. Suppose that $1\leq b_{\vphi}<+\infty$. Then, for any $\delta>0$ we have
\[
\liminf_{n\to\infty}
\frac{\log\log \vphi_{e,b_{\vphi} + \delta} (n)}{n}
=
\log(b_{\vphi}+\delta)
>
\log b_{\vphi}
=
\liminf_{n\to\infty}
\frac{\log\log \vphi(n)}{n}
\]
and, thus, $\clE_{\ell}''(\vphi_{e,b_{\vphi} + \delta}) \subseteq \clE_{\ell}'(\vphi)$. 
Then, Theorem \ref{TEO:luczak for primes} gives
\[
\frac{1}{1+b_{\vphi} + \delta}
\leq 
\hdim \clE_{\ell}'(\vphi),
\:\text{ so }\;
\frac{1}{1+b_{\vphi}} 
\leq \hdim \clE_{\ell}'(\vphi).
\]
When $b_{\vphi}=\infty$, the result follows from $\clE_{\ell}'(\vphi)\subseteq \clE_{\ell}(\vphi)$ and Theorem \ref{TEO:HWX:HDIM}.
\section{Further research}{\label{Section:FurtherResearch}}
We close the paper with a list of open problems for further research.
\subsection{Laws of large numbers}
For any $\ell,n\in\N$ and $x\in (0,1)$ irrational, define
\[
S_{n,\ell}(x)
=
\sum_{j=1}^n a_j(x)\cdots a_{j+\ell -1}(x)
\;\text{ and }\;
M_{n,\ell}(x) = \max_{1\leq j\leq n} a_j(x).
\]
It is well known that the partial quotients are not stochastically independent (with respect to the Gauss measure) nor integrable, so the classical laws of large numbers do not apply. 
There are, however, several results in this direction.
Khinchin proved a weak law of large numbers for the sequence of partial quotients \cite{Khintchine1935}. 
\begin{teo01}[{\cite[p. 377]{Khintchine1935}}]\label{THM:Khinchin:WLLN}
The sequence $\left( \frac{S_{n,1}}{n\log n}\right)_{n\geq 1}$ converges in Lebesgue measure to $\frac{1}{\log 2}$; that is, for every $\veps>0$ we have
\[
\lim_{n\to\infty}
\leb\left(\left\{x\in(0, 1): \left|\frac{ S_{n,1}(x)}{n\log n}-\frac1{\log 2}\right|\geq \epsilon\right\}\right)
=
0.
\]
\end{teo01}
The Borel-Bernstein Theorem tells us that almost every $x\in (0,1)$ satisfies $a_n(x)\geq n \log^2 n$ infinitely often, so 
\[
\limsup_{n\to\infty}
\frac{S_{n,1}(x)}{n\log n}
=
\infty
\;\text{ almost everywhere}.
\]
In fact, under some reasonable regularity properties, Philipp \cite[Theorem 1]{Philipp88} proved that there is no function $\sigma:\N\to\RE_{>0}$ such that $S_n(x)/\sigma(n)$ converges almost everywhere to a finite nonzero constant.
Diamond and Vaaler \cite{DiamondVaaler} showed that the convergence almost everywhere occurs after removing the largest term.
\begin{teo01}[{\cite[Corollary 1]{DiamondVaaler}}]\label{THM:DV:SLLN}
For almost every $x\in (0,1)$, we have
\[
\lim_{n\to\infty}
\frac{S_{n,1}(x) - M_{n,1}(x)}{n\log n}
=
\frac{1}{\log 2}.
\]
\end{teo01}
Hu, Hussain, and Yu \cite{HHY} obtained analogues of Theorems~\ref{THM:Khinchin:WLLN} and ~\ref{THM:DV:SLLN} for $\ell=2$.
\begin{teo01}\label{THM:HHY:LLN}\, \, 
\begin{enumerate}[i.]
\item \cite[Theorem 1.4]{HHY} The sequence $\left( \frac{S_{n,2}}{n\log^2 n}\right)_{n\geq 1}$ converges in Lebesgue measure to $\frac{1}{2 \log 2}$.
\item \cite[Theorem 1.5]{HHY} For almost every $x\in (0,1)$, we have 
\[
\lim_{n\to\infty}
\frac{S_{n,2}(x) - M_{n,2}(x)}{n\log^2 n}
=
\frac{1}{2\log 2}.
\]
\end{enumerate}
\end{teo01}
Among other probabilistic results, Schindler and Zweim{\"u}ller obtained prime analogues of Theorem \ref{THM:Khinchin:WLLN} and \ref{THM:HHY:LLN}. To state them, for $n,\ell\in\N$ and $x\in (0,1)\setminus \QU$ define
\[
S_{n,\ell}'(x)
=
\sum_{j=1}^n a_j'(x)\cdots a_{j+\ell -1}'(x)
\;\text{ and }\;
M_{n,\ell}'(x) = \max_{1\leq j\leq n} a_j'(x).
\]
\begin{teo01}\label{THM:ScZw:LLN}
\,\,
\begin{enumerate}[i.]
\item \label{THM:ScZw:LLN:i} \cite[Theorem 3.1]{SchindlerZweimuller2023} The sequence $\left( \frac{S_{n,1}'}{n\log \log n}\right)_{n\geq 1}$ converges in Lebesgue measure to $\log 2$.
\item \label{THM:ScZw:LLN:ii} \cite[Theorem 2.2]{SchindlerZweimuller2023} For almost every $x\in (0,1)$, we have 
\[
\lim_{n\to\infty}
\frac{S_{n,1}'(x) - M_{n,1}'(x)}{n\log\log n}
=
\frac{1}{\log 2}.
\]
\end{enumerate}
\end{teo01}
\begin{rema}
Theorem 3.1 in \cite{SchindlerZweimuller2023} is stronger than Theorem \ref{THM:ScZw:LLN}.\ref{THM:ScZw:LLN:i}. 
To deduce our statement from it, simply recall that convergence in measure is equivalent to convergence in distribution when the limit is constant \cite[Proposition 8.5.2]{MR3135152}.
\end{rema}
These developments motivate the following problems.
\begin{prob01}
Determine a function $\sigma:\N\to\RE_{>0}$ such that $\left( \frac{S_{n,2}'}{\sigma(n)}\right)_{n\geq 1}$ converges in measure to a positive finite constant and, for some positive finite constant $\kappa$,
\[
\lim_{n\to\infty}
\frac{S_{n,2}'(x) - M_{n,2}'(x)}{\sigma(n)}
=
\kappa
\quad\text{ almost everywhere}.
\]
\end{prob01}
Naturally, it would be desirable to extend Theorem~\ref{THM:Khinchin:WLLN} and \ref{THM:DV:SLLN} (or Theorem \ref{THM:ScZw:LLN}) to arbitrary $\ell\in \N$.
This problem, however, poses significant technical difficulties even without the primality restriction. 
A key step in the current proofs of Theorem~\ref{THM:Khinchin:WLLN}-\ref{THM:HHY:LLN} is calculating, for a particular function $\vphi$, the Lebesgue measure of the set of irrational numbers $x\in (0,1)$ such that for infinitely many $n\in\N$ there are integers $i,j$ satisfying
\[
1\leq i< j \leq n, 
\quad
\prod_{r=0}^{\ell-1} a_{j+r}(x)\geq \vphi(n),
\quad \text{ and }\quad
\prod_{r=0}^{\ell-1} a_{i+r}(x)\geq \vphi(n).
\]
The problem for non-decreasing $\vphi$ satisfying $\vphi(n)\to \infty$ as $n\to\infty$ is solved for $\ell=1$ \cite{MR4577484}, $\ell=2$ \cite{MR4687014}, and $\ell=3$ \cite{arXiv240510538}.
\begin{prob01}
Let $\vphi:\N\to\mathbb{R}_{>0}$ be a non-decreasing function such that $\vphi(n)\to\infty$ as $n\to\infty$.
Compute the Lebesgue measure and the Hausdorff dimension of the set of irrational numbers $x\in (0,1)$ such that for infinitely many $n\in\N$ we have
\[
1\leq i< j \leq n, 
\quad
\prod_{r=0}^{\ell-1} a_{j+r}'(x)\geq \vphi(n),
\quad \text{ and }\quad
\prod_{r=0}^{\ell-1} a_{i+r}'(x)\geq \vphi(n).
\]
\end{prob01}
The statements \cite[Theorem 1.6-1.8]{HHY} provide a finer description of the sets of irrational numbers where $(S_{n,2}(x))_{n\geq 1}$ grows at a prescribed rate. The corresponding problem in the prime setting reads as follows:
\begin{prob01}
For any non-decreasing function $\vphi:\N\to\RE_{>0}$ and $\ell\in\N$, determine the Hausdorff dimension of the set 
\[
E_{\ell}(\varphi)
:=
\left\{x\in(0,1):\lim\limits_{n\to\infty}\frac{S_{n,\ell}'(x)}{\varphi(n)}=1\right\}.
\]
\end{prob01}
\subsection{Miscellaneous problems}
First, we may replace the set of prime numbers $\clP$ with any infinite subset $\clA$ of natural numbers.
\begin{prob01}
For any $\ell\in\N$, any positive function $\varphi:\N\to\RE_{>0}$ and any infinite subset $\clA\subseteq \N$, define 
\[
\clE_{\ell}(\vphi;\clA)
=
\left\{
x\in (0,1):
a_n(x)\cdots a_{n+\ell-1}(x)\geq \vphi(n)
\;\text{ and }\;
a_n(x)\cdots a_{n+\ell-1}(x)\in\clA
\;\text{ for i.m.} n\in\N
\right\}.
\]
\end{prob01}
\begin{rema}
The proofs of Theorem~\ref{THM:LEBMEAS} and \ref{TEO:MAIN:HDIM} rely on the Prime Number Theorem. 
Hence, we might need some regularity conditions on the function $\N\to\N_0$, $n\mapsto \#\{1,2,\ldots, n\}\cap \clA$.
\end{rema}
Khinchin proved there is a constant $K= 2.68545^{+}$ such that for almost every $x\in (0,1)$ we have
\[
\lim_{n\to\infty} \left(a_1(x)\cdots a_n(x)\right)^{\frac{1}{n}}
=
K.
\]
For a proof, see, for example, \cite[Corollary 3.8]{EinsiedlerWardBook}.
\begin{prob01}
For every $\alpha\in (0,\infty]$, determine the Hausdorff dimension of the set 
\[
\clK'(\alpha)
:=
\left\{ 
x\in (0,1): \lim_{n\to\infty} \left(a_1'(x)\cdots a_n'(x)\right)^{\frac{1}{n}}=\alpha
\right\}.
\]
\end{prob01}
\begin{rema}
We assume $\alpha>0$ because $\left(a_1'(x)\cdots a_n'(x)\right)^{\frac{1}{n}}=0$ for large $n$ almost everywhere. This problem was addressed without the primality assumption in \cite{MR2470627}.    
\end{rema}
Let $\Phi:\N\to\RE_{>0}$ be any positive function such that $\Phi(n)\to\infty$ as $n\to\infty$.
One way of understanding the behaviour of the product of prime partial quotients is to investigate the Hausdorff dimension of the following sets:
\begin{equation*}
    \Lambda_{\rm inf}(\Phi):=\left\{x\in[0, 1): \liminf_{n\to\infty}\frac{\log\left(1+\displaystyle\prod_{i=1}^ma^\prime_{n+i}(x)\right)}{\Phi(n)}=1\right\},
\end{equation*}
and
\begin{equation*}
    \Lambda_{\rm sup}(\Phi):=\left\{x\in[0, 1): \limsup_{n\to\infty}\frac{\log\left(1+\displaystyle\prod_{i=1}^ma^\prime_{n+i}(x)\right)}{\Phi(n)}=1\right\}.
\end{equation*}
Clearly, $\leb(\Lambda_{\rm inf}(\Phi))=0$ by Birkhoff's ergodic theorem.
Similar to the above sets, for any $0\leq \alpha\leq \beta\leq \infty$ define
\begin{equation*}
    \Lambda_{\alpha, \beta}:=\left\{x\in[0, 1): \liminf_{n\to\infty}\frac{\log\left(1+\displaystyle\prod_{i=1}^ma^\prime_{n+i}(x)\right)}{\log q_n(x)}=\alpha,  \limsup_{n\to\infty}\frac{\log\left(1+\displaystyle\prod_{i=1}^ma^\prime_{n+i}(x)\right)}{\log q_n(x)}=\beta\ \right\}.
\end{equation*}
\begin{prob01}
For $\Phi$, $\alpha$, and $\beta$ as above, compute the Hausdorff dimension of $\Lambda_{\rm inf}(\Phi)$, 
$\Lambda_{\rm sup}(\Phi)$,  and $\Lambda_{\alpha, \beta}$.
\end{prob01}
\begin{rema}
In the previous sets, we add 1 to avoid evaluating $\log$ at $0$ whenever the partial quotients are composite.
Moreover, assuming every partial quotient is prime, the limits are unaffected by adding $1$.
\end{rema}

%
\end{document}